\documentclass[11pt,a4paper]{article}

\usepackage{amssymb}
\usepackage{latexsym}
\usepackage{amsfonts}
\usepackage{mathrsfs}
\usepackage{amscd}
\usepackage{amsmath}
\usepackage{amssymb}
\usepackage{wasysym}
\usepackage{bm}
\newtheorem{Th}{Theorem}
\newtheorem{Prop}{Proposition}

\newtheorem{Lma}{Lemma}[section]
\newtheorem{Dfi}{Definition}
\newtheorem{Rm}{Remark}

\newcommand{\be}{\begin{equation}}
\newcommand{\ee}{\end{equation}}
\newcommand{\bes}{\begin{equation*}}
\newcommand{\ees}{\end{equation*}}
\newcommand{\R}{\mathbb{R}}

\newcommand{\C}{\mathbb{C}}

\newcommand\res{\mathop{\hbox{\vrule height 7pt width .5pt depth 0pt
\vrule height .5pt width 6pt depth 0pt}}\nolimits}

\newcommand{\reset}{\setcounter{equation}{0}\setcounter{Th}{0}\setcounter{Prop}{0}\setcounter{Co}{0}\setcounter{Lm}{0}\setcounter{Rm}{0}}

\def\ti{\tilde}
\def\lf{\left}
\def\rg{\right}

\def\al{\alpha}
\def\la{\lambda}
\def\eps{\varepsilon}

\def\Om{\Omega}

\def\p{\partial}
\def\pro{\pi_{\vec{n}}}
\def\bn{\vec{n}}
\def\bna{\vec{n}_\al}

\def\bex{\bbe_1}
\def\bey{\bbe_2}
\def\bez{\bbe_z}
\def\bea{\bbe_a}
\def\beb{\bbe_b}
\def\bei{\bbe_i}
\def\bej{\bbe_j}

\def\bezz{\bbe_{z^*}}
\def\px{\partial_{x_1}}
\def\py{\partial_{x_2}}
\def\pj{\partial_{x_j}}
\def\pk{\partial_{x_k}}
\def\pa{\partial_{a}}
\def\pb{\partial_{b}}
\def\pz{\partial_{z}}
\def\pzz{\partial_{z^*}}
\def\bbe{\vec{e}}
\def\bH{\vec{H}}

\def\bv{\vec{v}}

\def\ba{\vec{a}}
\def\bbe{\vec{e}}

\def\bL{\vec{L}}
\def\bR{\vec{R}}
\def\bV{\vec{V}}

\def\bX{\vec{X}}
\def\bY{\vec{Y}}
\def\bw{\vec{w}}

\def\bB{\vec{B}}

\def\bp{\vec{\Phi}}
\def\bP{\vec{\Phi}}
\def\bxi{\vec{\xi}}

\def\nablap{\nabla^\perp}
\def\bet{\beta}
\def\bul{\bullet}
\def\di{\mathbb{D}^2}
\def\res{\mathop{\hbox{\vrule height 7pt width .5pt 
depth 0pt\vrule height .5pt width 6pt depth 0pt}}\nolimits}

\begin{document}

\reset

\title{Local Palais-Smale Sequences\\
for the Willmore Functional}
\author{Yann Bernard$^{*}$\:,\:Tristan Rivi\`ere\footnote{Department of Mathematics, ETH Zentrum,
8093 Z\"urich, Switzerland.}}
\maketitle
$\textbf{Abstract :}$ 
Using the reformulation in divergence form of the Euler-Lagrange equation for the Willmore functional as it was developed in the second author's paper \cite{Riv2}, we study the limit of a local Palais-Smale sequence of weak Willmore immersions with locally square-integrable second fundamental form. We show that the limit immersion is smooth and that it satisfies the {\it conformal Willmore equation}: it is a critical point of the Willmore functional restricted to infinitesimal conformal variations.  
\section{Introduction}

\subsection{The Willmore Functional and the Willmore Equation}

Consider an oriented surface $\Sigma$ without boundary immersed in $\mathbb{R}^m$, for some $m\ge3$, through the action of a smooth positive immersion $\bP$. We introduce the Gau\ss\, map $\bn$, which to every point $x$ in $\Sigma$ assigns the unit $(m-2)$-plane $N_{\bP(x)}\bP(\Sigma)$ orthogonal to the oriented tangent space $T_{\bP(x)}\bP(\Sigma)$. This map acts from $\Sigma$ into $Gr_{m-2}(\mathbb{R}^m)$, the Grassmannian of oriented $(m-2)$-planes in $\mathbb{R}^m$. It thus naturally induces a projection map $\pi_{\bn}$ which to every vector $\xi$ in $T_{\bP(x)}\mathbb{R}^m$ associates its orthogonal projection $\pi_{\bn}(\xi)$ onto $N_{\bP(x)}\bP(\Sigma)$. Let $x$ be a point on $\Sigma$. We denote by $\bB_x$ the second fundamental form of the immersion $\bP$. It is a $N_{\bP(x)}\bP(\Sigma)$-valued symmetric bilinear form on $T_x\Sigma\times T_x\Sigma$, defined by $\bB_x=\pi_{\bn}\circ d^2\bP$. Having chosen an orthonormal basis $\{e_1,e_2\}$ on $T_x\Sigma$, the mean curvature $\bH(x)$ of the immersion $\bP$ is the vector in $N_{\bP(x)}\bP(\Sigma)$ given by
\bes
\bH(x)\;:=\;\dfrac{1}{2}\,\text{Tr}\,(\bB_x)\;\equiv\;\dfrac{1}{2}\,\big(\bB_x(e_1,e_1)\,+\,\bB_x(e_2,e_2)\big)\:.
\ees
The {\it Willmore functional} is the Lagrangian
\be\label{wil}
W\big(\bP(\Sigma)\big)\;:=\;\int_{\Sigma}\,\big|\bH\big|^2\,d\mu_g\:,
\ee
where $d\mu_g$ is the area form of the metric $g$ induced on $\bP(\Sigma)$ via the canonical metric on $\mathbb{R}^m$.\\[2ex]
The critical points of (\ref{wil}) for perturbations of the form $\bP+t\,\bxi$\,, where $\bxi$ is an arbitrary compactly supported smooth map on $\Sigma$ into ${\R}^m$ are known as {\it Willmore surfaces}. Examples of Willmore surfaces are legions. Any minimal surface, i.e. one for which $\bH\equiv0$\,, realizes an absolute minimum of the Willmore functional. Round spheres are also Willmore surfaces, and, more generally, all Willmore surfaces with genus zero were obtained by Robert Bryant \cite{Bry1}, \cite{Bry2}. Another important example was devised by Willmore  in 1965. It is the torus of revolution obtained through rotating a circle of radius 1 whose center is located at a distance $\sqrt{2}$ from its axis of rotation (equivalently, it is the stereographic projection into $\mathbb{R}^3$ of the Clifford torus). Willmore proved that this torus is indeed a Willmore surface, and he conjectured that it minimizes, up to M\"obius transforms, the Willmore energy in the class of smooth and immersed tori. Despite partial answers, this assertion (known as {\it the Willmore conjecture}) remains unsolved. For more details, the reader is referred to \cite{BK}, \cite{LY}, \cite{Sim}, and \cite{Wil2}. Further examples of Willmore surfaces are profuse in the literature, and we content ourselves with citing \cite{Wil2}, \cite{PS}, and the references therein.\\[1ex]
The Euler-Lagrange equation obtained through varying the Willmore functional as aforementioned was first\footnote{although it seemingly was known to Thomsen and Blaschke decades sooner.} derived in \cite{Wil1} in the three-dimensional case, and subsequently extended by Weiner \cite{Wei} for general $m\ge3$. We now recall this equation.\\
Given any vector $\bw$ in $N_{\bP(x)}\bP(\Sigma)$\,, consider
the symmetric endomorphism $A_x^{\bw}$ of $T_x\Sigma$ satisfying $g(A_x^{\bw}(\bX),\bY)=B_x(\bX,\bY)\cdot \bw$\,, where $\cdot$ denotes the standard scalar product in ${\R}^m$\,, for every pair of vectors $\bX$ and $\bY$ in $T_x\Sigma$. 
The map $A_x:\bw\mapsto A^{\bw}_x$\, is a homomorphism from $N_{\bP(x)}\bP(\Sigma)$ into the linear space of symmetric endomorphisms on $T_x\Sigma$. We next define $\ti{A}_x=$ $ ^{t}A_x\circ A_x$\,, which is an endomorphism of $N_{\bP(x)}\bP(\Sigma)$.
If $\{e_1,e_2\}$ is an orthonormal basis
of $T_x\Sigma$\,, and if $\bL$ is a vector in $N_{\bP(x)}\bP(\Sigma)$\,, then it is readily seen that
$\ti{A}(\bL)=\sum_{i,j} \bB_x(e_i,e_j)\bB_x(e_i,e_j)\cdot\bL$.
With this notation, as shown in \cite{Wei}, $\bP$ is  a smooth Willmore immersion if and only if it satisfies Euler-Lagrange equation 
\be\label{wileq} 
\Delta_\perp\bH\,-\,2\,|\bH|^2\bH\,+\,\ti{A}(\bH)\:=\:0\:,
\ee
where $\Delta_\perp$ is the negative covariant Laplacian for the connection $D$ in the normal bundle $N\bP(\Sigma)$\, derived from the ambient scalar product in ${\R}^m$. Namely, for every section $\sigma$ of $N\bP(\Sigma)$\,, one has $D_{\bX}\sigma:=\pi_{\bn}(\sigma_\ast\bX)$. \\

The paternity of the Willmore functional is delicate to establish precisely. Although it bears the name of T.J. Willmore whom studied it in 1965 \cite{Wil1} thereby initiating its popularization, the Willmore functional had been previously considered in the works of Sophie Germain \cite{Ger}, Gerhard Thomsen \cite{Tho}, and Wilhelm Blaschke \cite{Bla}.  Its wide range of applications includes various areas of science, where it plays an important r\^ole. Amongst others, the Willmore functional appears in molecular biophysics as the surface energy for lipid bilayers in the Helfrich model \cite{Hef} (cell membranes tend to position themselves so as to minimize the Willmore energy) ; in solid mechanics as the limit-energy for thin plate theory \cite{FJM} ; in general relativity as the main contributing term to the Hawking quasilocal mass (cf. \cite{Haw}, \cite{HI}) ; in string theory as an extrinsic string action \`a la Polyakov \cite{Pol}. \\[1.5ex]
The importance of the Willmore functional is largely due to its invariance under conformal transformations of the metric of the ambient space. This remarkable property was first brought into light by White \cite{Whi} in the three-dimensional case, then generalized by B.Y. Chen \cite{Che}. As the second author of the present paper showed in \cite{Riv1}, Euler-Lagrange equations arising from a two-dimensional conformally invariant Lagrangian with quadratic growth can be written in divergence form. These ``conservations laws" fostered within variational problems involving conformally invariant Lagrangians offer a significant help. In particular, the general ideas introduced in \cite{Riv1} are led to fruition in \cite{Riv2}, where conservation laws relative to the Willmore functional are developed and successfully applied to produce a variety of interesting results. Our present work stems from this alternative formulation of the Willmore equation. 

\subsection{Weak Willmore Immersions and Conservation Laws}

This section is devoted to recalling the formalism introduced in \cite{Riv2} and the results therein established, which compose the foundation of our work. \\[-1ex]

Owing to the Gau\ss-Bonnet theorem, we note that the Willmore energy (\ref{wil}) may be equivalently expressed as
\bes
W\big(\bP(\Sigma)\big)\,=\,\int_{\Sigma}\,\big|\bB\big|^2\,d\mu_g\;+\;4\,\pi\,\chi(\Sigma)\:,
\ees
where $\bB$ is the second fundamental form, and $\chi(\Sigma)$ is the Euler characteristic of $\Sigma$. Since the latter is a topological invariant, from the variational point of view, we infer that Willmore surfaces are the critical points of the energy
\bes
\int_{\Sigma}\,\big|\bB\big|^2\,d\mu_g\:.
\ees
In studying such surfaces, it thus appears natural to restrict our attention on immersions whose second fundamental forms are locally square-integrable. More precisely, we work within the framework of {\it weak immersion with locally $L^2$-bounded second fundamental form}, which we now define.

\smallskip
\begin{Dfi}
Let $\bP$ be a $W^{1,2}$-map from a two-dimensional manifold $\Sigma$ into ${\R}^m$.
$\bP$ is called a weak immersion with locally $L^2$-bounded second fundamental form
whenever there exist locally about every point an open disk $D$, a positive constant $C$, and a sequence of smooth embeddings $(\bP_k)$ from $D$ into ${\R}^m$, such that
\begin{itemize}
\item[i)] $\hspace{4cm}\displaystyle{\mathcal{H}_2(\bP(D))\,\ne\,0}$\:,
\item[ii)] $\hspace{4cm}\displaystyle{{\mathcal H}_2(\bP_k(D))\,\le\, C\,<\,\infty}$\:,
\item[iii)] $\hspace{4cm}\displaystyle{\int_{D}|B_k|^2\ dvol_{g_k}\,\le\,\frac{8\pi}{3}}$\:,
\item[iv)]$\hspace{4cm}\displaystyle{\bP_k\rightharpoonup\bP\quad\mbox{weakly in}\:\:W^{1,2}}$\:,\\[-2ex]
\end{itemize}
where ${\mathcal H}_2$ is the two-dimensional Hausdorff measure, $B_k$
is the second fundamental form associated to the embedding $\bP_k$, and $g_k$ denotes the metric on $\bP_k(\Sigma)$ obtained via the pull-back by $\bP_k$ of the induced metric.
\end{Dfi}

\noindent
A useful characterization of weak immersions with square-integrable second fundamental form was originally obtained by Tatiana Toro \cite{To1} and \cite{To2}, and by Stefan M\"uller and Vladim\'ir \v{S}ver\`ak in \cite{MS}. In the same spirit, Fr\'ed\'eric H\'elein obtained the following statement (whose proof appears as that of Theorem 5.1.1 in \cite{Hel}). 

\begin{Th}
\label{th-I.2a} 
Let $\bP$ be a weak immersion from a two-dimensional manifold $\Sigma$ into ${\R}^m$ with locally $L^2$-bounded second fundamental form. Then locally about every point on $\Sigma$, there exist an open disk $D$ and a homeomorphism $\Psi$ of $D$ such that $\bP\circ\Psi$ is a conformal bilipschitz immersion. In this parametrization, the metric $g$ on $D$ induced by the standard metric of ${\R}^m$ is continuous. Moreover, the Gau\ss\, map $\bn$ of this immersion lies in $W^{1,2}(D,Gr_{m-2}({\R}^m))$, relative to the induced metric $g$.
\end{Th}

We previously observed that weak immersions with square-integrable second fundamental form are particularly suited for the study of Willmore surfaces. More importantly, in \cite{Riv2}, the notion of {\it weak Willmore immersion} is introduced. It is based on the following theorem, established in the same paper. It provides a reformulation of the Willmore equation in divergence form.  

\begin{Th}
\label{th-0.1}
The Willmore equation (\ref{wileq}) is equivalent to 
\bes
d\lf(\ast_g\, d\bH-3\ast_g\pi_{\bn}\big(d\bH\big)\rg)\,-\,d\star\lf( d\bn\wedge \bH\rg)\:=\:0\:,
\ees
 where $\ast_g$ is the Hodge operator on $\Sigma$ associated with the induced metric $g$, and $\,\star$ is the usual Hodge operator on forms. \\[1ex]
\noindent
In particular, a conformal immersion $\bP$ from the flat disc $\di$ into ${\R}^m$ is Willmore if and only if
\be\label{conslaw}
\Delta\bH\,-\,3\ div\big(\pi_{\bn}(\nabla \bH)\big)\,+\,div\star\lf(\nabla^\perp\bn\wedge\bH\rg)\:=\:0\:,
\ee
where the operators $\nabla$, $\nabla^\perp$, $\Delta$, and $div$ are understood with respect to the flat metric on $\di$. Namely, $\nabla=(\p_{x_1},\p_{x_2})$, $\nabla^\perp=(-\p_{x_2},\p_{x_1})$, $\Delta=\nabla\cdot\nabla$, and \,$div =tr\circ\nabla$.
\end{Th}

\noindent
We are now ready to define the notion of {\it weak Willmore immersion with $L^2$-bounded second fundamental form}.

\begin{Dfi}
\label{df-I.3a}
A  weak immersion $\bP$ from a two-dimensional manifold $\Sigma$ into ${\R}^m$ with 
locally $L^2$-bounded second fundamental form is Willmore whenever (\ref{conslaw}) holds in the sense of distributions locally about every point in a conformal parametrization from the two-dimensional disk $D$, as indicated in Theorem~\ref{th-I.2a}.
\end{Dfi}

Henceforth, we will work with weak Willmore immersions, which we shall assume to be conformal on the unit-disk. 

\bigskip

We introduce the local coordinates $(x_1,x_2)$ for the flat metric on the unit-disk $\,\di=\big\{(x^1,x^2)\in{\R}^2\ ;\ x_1^2+x_2^2<1\big\}$. The operators $\nabla$, $\nabla^{\perp}$, $div$, and $\Delta$  will be understood in these coordinates.
Let $\bp:\di\rightarrow\R^m$ be a conformal immersion. We define the conformal factor $\lambda$ via
\bes
\px\bp\:=\:\text{e}^{\la}\:=\:\py\bp\:.
\ees
Thanks to the topology of $\di$, there exists, for almost every $z\in\di$, a positively oriented orthonormal basis $\{\bn_1,\ldots,\bn_{m-2}\}$ of $N_{\bp(z)}\bp(\di)$, complement of the tangent plane to $\bp(\di)$ at $\bp(z)$, such that\\
$\{\bex,\bey,\bn_1,\ldots,\bn_{m-2}\}$ forms a basis of $T_{\bp(z)}\R^m$. Following the ``Coulomb gauge extraction method" exposed in the proof of Lemma 4.1.3 from \cite{Hel}, the basis $\{\bn_\al\}_{\al=1}^{m-2}$ may be chosen to satisfy 
\bes
div\,\big\langle\nabla\bn_\al\,,\bn_\beta\big\rangle\:=\:0\:,\qquad\forall\:\:\:1\,\le\,\al\,,\beta\,\le\,m-2\:.
\ees
From the Pl\"ucker embedding, which realizes $Gr_{m-2}(\R^m)$ as a submanifold of the projective space of the $(m-2)^\text{th}$ exterior power $\,\mathbb{P}\big(\bigwedge^{m-2}\R^m\big)$, we can represent the Gau\ss\, map as the $(m-2)$-vector
\bes
\bn\:=\:\bigwedge_{\al=1}^{m-2}\,\bn_\al\:.
\ees
Via the Hodge operator $\star\,$, we may identify vectors and $(m-1)$-vectors in ${\R}^m$. In particular, there holds
\bes
\star\,(\bn\wedge \bbe_1)\;=\;\bbe_2\qquad\mbox{and}\qquad\star(\bn\wedge \bbe_2)=-\,\bbe_1\:.
\ees
With this notation, the second fundamental form $\bB$, which is a symmetric 2-form on $T_{\bp(z)}\bp(\di)$  into  $N_{\bp(z)}\bp(\di)$, may be expressed as
\be\label{vecb}
\bB\;=\;\sum_{\al,i,j}\ \text{e}^{-2\la}\,h^\al_{ij}\ \bn_\al\,dx^i\otimes dx^j\;\equiv\;\sum_{\al,i,j}\ h^\al_{ij}\ \bn_\al\,(\bbe_i)^\ast\otimes(\bbe_j)^\ast\:,
\ee
with
\bes
h^\al_{ij}\;=\;-\,\text{e}^{-\la}\,\bei\cdot\pj\bn_\al\:.
\ees\\
The mean curvature vector $\bH$ is
\bes
\label{II.2}
\bH\:=\:\sum_{\al=1}^{m-2}\,H^\al\,\bn_\al\:=\:\;\frac{1}{2}\,\sum_{\al=1}^{m-2}\,\big(h^\al_{11}+h^\al_{22}\big)\, \bn_\al\:.
\ees
The Weingarten operator $\bH_0$ is
\bes
\label{II.20}
\bH_0\:=\:\sum_{\al=1}^{m-2}\,H_0^\al\,\bn_\al\:=\:\;\frac{1}{2}\,\sum_{\al=1}^{m-2}\,\big(h^\al_{11}-h^\al_{22}\,+\,2\,i\,h^\al_{12}\big)\, \bn_\al\:.
\ees
\medskip
\noindent
In this framework, the Euler-Lagrange equation (\ref{wileq}) for the Willmore functional is cast in the form 
\be\label{wil1}
\Delta_\perp\bH\,+\,\sum_{\al,\beta,i,j}\,h^\al_{ij}\,h^\bet_{ij}\,H^\bet\,\bn_\al\:-\:2\,\big|\vec{H}\big|^2\bH\:=\:0\:,
\ee
with
\bes
\Delta_\perp\bH\:=\:\text{e}^{-2\la}\,\pi_{\bn}\,div\big(\pi_{\bn}(\nabla\bH)\big)\:.
\ees

\reset

\section{Main Results}

\subsection{Local Palais-Smale Sequences}

As we stated in the Introduction, T. Willmore conjectured in 1965 that the Willmore torus (defined on page 2) minimizes, up to M\"obius transformations, the Willmore energy in the class of smooth immersed tori in $\R^3$. To this day, no satisfactory demonstration of this assertion has been found. Amid the works aimed at solving this problem, one fundamental property of the Willmore functional was brought into light by Leon Simon \cite{Sim}. Namely, for each dimension $m\ge3$, there exists a compact embedded real analytic torus in $\mathbb{R}^m$ which minimizes the Willmore energy in the class of compact, genus 1, embedded surfaces without boundary\footnote{more generally, Simon obtains an analogous statement for each genus $\frak{g}\in\mathbb{N}$.}. Unfortunately, it remains unknown whether this minimizer is the Willmore torus. Simon's rather sophisticated proof is constructive in nature: an explicit minimizing sequence for the Willmore functional is devised. This is thus one instance in which studying and understanding minimizing sequences of the Willmore functional are essential. Yet more generally, it is interesting to investigate Palais-Smale sequences for the Willmore functional. Such will be our goal in this section. \\

We open our considerations with an ``empirical" observation which will hopefully convince the reader that the results derived in \cite{Riv2} offer a suitable framework to acquire information on Palais-Smale sequences of Willmore surfaces. \\
Let us consider a conformal weak Willmore immersion $\bp\,$ from the flat disk $\di$ into $\mathbb{R}^m$ with bounded square-integrable second fundamental form. Up to an additive constant, it is possible to define a map $\bL$ satisfying\footnote{refer to the Introduction or the Appendix for the notation.}
\be\label{defL}
\nabla^\perp\bL\::=\:\nabla^\perp\bH\,-\,3\,\pi_{\bn}\big(\nabla\bH\big)\,+\,\star\,\big(\nabla^\perp\bn\wedge\bH\big)\:.
\ee
In \cite{Riv2}, it is shown that the following system holds:
\be\label{wil2bis}
\left\{\begin{array}{rcl}
\nabla\bp\,\cdot\nabla^\perp\bL&=&0\\[1ex]
\nabla\bp\,\wedge\nabla^\perp\bL&=&2\,(-1)^m\;\nabla\big(\star(\bn\,\res\bH)\big)\,\res\nabla^\perp\bp\:.
\end{array}\right.
\ee

\noindent
The Hodge operator $\star$ and the contraction operator $\res$ commute with the partial differentiation operators. Accordingly, the terms appearing in the system (\ref{wil2bis}) enjoy a peculiar property: more than mere ``products of derivatives", they can be factored in divergence form. This structural feature has the analytical advantage of being robust under weak limiting process. Hence it is legitimate to hope that a local Palais-Smale sequence of conformal weak Willmore immersions with, say, uniformly bounded square-integrable second fundamental forms, converges to an element $\bp$ satisfying the system (\ref{wil2bis}) for some function $\bL$ (which may or may not be related to $\bp$ via (\ref{defL})). \\[1ex]
This is essentially the result that we shall establish. Prior to stating precisely, it is necessary to define the notion of local Palais-Smale sequence for the Willmore functional.  

\begin{Dfi}
Let $(\bp_k)$ be a sequence of conformal immersions from the unit-disk $\di$ into $\mathbb{R}^m$ such that
\be\label{ps1}
\big\Vert\bp_k\big\Vert_{W^{2,2}\,\cap\,W^{1,\infty}}\:\le\:C
\ee
holds uniformly for some positive constant $C$. Denoting respectively by $\bn_k$ and $\bH_k$ the Gau\ss\, map and the mean curvature vector associated with the immersion $\bp_k$, we set
\bes
Q_k\;:=\;\nabla\bH_k\,-\,3\,\pi_{\bn_k}\big(\nabla\bH_k\big)\,+\,\star\,\big(\nabla^\perp\bn_k\wedge\bH_k\big)\:.
\ees
The sequence $(\bp_k)$ is locally Palais-Smale if, in addition to (\ref{ps1}), it satisfies
\bes
div\,Q_k\:\longrightarrow\:0\qquad\text{strongly\: in}\quad \big(W^{2,2}\cap W^{1,\infty}\big)'(\di)\:.
\ees
\end{Dfi}

\medskip
\noindent
The following result provides a first description of the limit of a Palais-Smale sequence of the Willmore functional.

\begin{Th}
\label{th-ps}
Let $(\bp_k)$ be a local Palais-Smale sequence of conformal immersions from the unit-disk\, $\di$ into $\mathbb{R}^m$. There exist a conformal weak immersion $\bp\in (W^{2,2}\cap W^{1,\infty})(\di)$ and an element $\bL\in L^{2,\infty}(\di)$ such that, up to extraction of a subsequence, 
\bes
\bp_k\:\longrightarrow\:\bp\qquad\text{in}\quad \mathcal{D}'(\di)
\ees
and the system
\be\label{wil2}
\left\{\begin{array}{rcl}
\nabla\bp\,\cdot\nabla^\perp\bL&=&0\\[1ex]
\nabla\bp\,\wedge\nabla^\perp\bL&=&-\,2\;\nabla\bp\,\wedge\nabla\bH
\end{array}\right.
\ee\\
holds in the sense of distributions, where $\bH$ denotes the mean curvature vector associated with $\bp$.
\end{Th}

\smallskip

\noindent
The apparent difference between the systems (\ref{wil2bis}) and (\ref{wil2}) is fictitious only. Indeed, one verifies\footnote{cf. equation (II.49) in the paper \cite{Riv2}.} the identity
\bes
\nabla\bp\,\wedge\nabla\bH\;=\;(-1)^{m-1}\;\nabla\big(\star(\bn\,\res\bH)\big)\,\res\nabla^\perp\bp\:.
\ees

\subsection{The Conformal Willmore Equation}

The convergence result stated in Theorem~\ref{th-ps} naturally begs the following question: if an immersion $\bp$ satisfies the system (\ref{wil2}), for some function $\bL$, is it true that $\bp$ is Willmore~? Unfortunately, and perhaps surprisingly, as far as the authors know, the answer is negative. Nevertheless, the following identification can be verified. 

\begin{Th}
\label{th-cw}
Let $\bp$ be a conformal weak immersion from the unit-disk $\di$ into $\R^m$ which lies in $\,W^{2,2}\cap W^{1,\infty}$ and such that
\bes
\int_{\di}\,\big|\nabla\bn\big|^2\:\le\:\eps\:,
\ees
for some $\eps>0\,$ small enough. 
Then there exists $\bL\in L^{2,\infty}(\di)$ such that $\bp$ satisfies (\ref{wil2}) if and only if $\bp$ is smooth and there holds
\be\label{cw}
\Delta_\perp\bH\,+\,\sum_{\al,\beta,i,j}\,h^\al_{ij}\,h^\bet_{ij}\,H^\bet\,\bn_\al\:-\:2\,\big|\vec{H}\big|^2\bH\:=\:e^{-2\la}\,\big\langle\vec{H}_0\,,f\big\rangle_{\mathbb{C}}\:,
\ee
for some holomorphic function $f$. 
\end{Th}

\begin{Rm}
\label{rem-cw}
An identity equivalent to (\ref{cw}) is obtained via setting
\be\label{defL0}
\nabla^\perp\bL_0\::=\:\nabla^\perp\bH\,-\,3\,\pi_{\bn}\big(\nabla\bH\big)\,+\,\star\,\big(\nabla^\perp\bn\wedge\bH\big)\:.
\ee
The statement of Theorem~\ref{th-cw} may then be rephrased through the equation
\be\label{cwbis}
\Delta\big(\bL-\bL_0\big)\:=\:2\,i\,\bH_0\,f\:.
\ee
\end{Rm}

The ``if" part of the statement of Theorem~\ref{th-cw} is clear. Indeed, choosing $f\equiv0$ and $\bL=\bL_0$ as in (\ref{defL0}), the validity of the system (\ref{wil2}) was settled in \cite{Riv2}. In the proof, we shall thus restrict our attention on the ``only if" part of the statement of Theorem~\ref{th-cw}. \\[-1ex]

Although the amount of information on the holomorphic function $f$ is rather limited, we note that the way it appears on the right-hand side of (\ref{cw}) makes it play the r\^ole of a Lagrange multiplier in the Willmore equation (\ref{wil1}). This observation enables one to give a natural geometric interpretation of $f$ in the context of Teichm\"uller theory (cf. \cite{BPP}).\\[0ex]

The Willmore functional being conformally invariant, the left-hand side of (\ref{wil1}), and thus of (\ref{cw}), remains unchanged under a holomorphic change of coordinates. An easy computation reveals that $\bH_0$ is the coefficient of $\,dz\otimes dz\,$ in the two-form $2\vec{B}$ given in (\ref{vecb}). Looking at the right-hand side of (\ref{cw}), we deduce  that $f$ must be the coordinate of a section $\,f(z)\,\pz\otimes\pz$\,. There is thus a one-to-one correspondence between the vector space to which $f$ belongs and the vector space $H^0(\Sigma)$ of holomorphic quadratic differentials. Let $\frak{g}$ be the genus of the immersed surface $\Sigma$. From the Riemann-Roch theorem (cf. Corollary 5.4.2 in \cite{Jo}), the space $H^0(\Sigma)$ satisfies
\bes
\dim_{\C}H^0(\Sigma)\;=\;\left\{\begin{array}{lcl}
0&,&\frak{g}=0\\
1&,&\frak{g}=1\\
3(\frak{g}-1)&,&\frak{g}\ge2\:.
\end{array}\right.
\ees

\medskip

The equation (\ref{cw}) is not a novelty ; it has been studied in various contexts. Immersions which satisfy (\ref{cw}) are sometimes known as {\it constrained Willmore immersions}, such as in \cite{BPP}. There, it is shown that (\ref{cw}) is the Euler-Lagrange equation deriving from the Willmore functional (\ref{wil}) under smooth compactly supported infinitesimal {\it conformal} variations. The corresponding critical points are the conformal-constrained Willmore surfaces. This notion clearly generalizes that of a Willmore surface, obtained via all smooth compactly supported infinitesimal variations. Conformal-constrained Willmore surfaces form a M\"obius invariant class of surfaces fostering remarkable properties, some of which are studied in \cite{BPP}, \cite{BPU}, \cite{Bry1}, and \cite{Ric}. In the latter, it is in particular established that every constant mean curvature surface in a 3-dimensional space-form is conformal-constrained Willmore. We have chosen to refer to (\ref{cw}) as the {\it conformal Willmore equation}, rather than as the somewhat vaguer term ``constrained Willmore equation". There are many ways to constrain the variations of the Willmore functional (e.g. restrictions on the volume and surface area in the Helfrich model). The adjective ``conformal" appears more descriptive in our situation.

\reset

\section{Proofs of the Theorems}

\subsection{Proof of Theorem~\ref{th-ps}}

Let us set
\bes
Q_k\;:=\;\nabla\bH_k\,-\,3\,\pro\big(\nabla\bH_k\big)\,+\,\star\,\big(\nabla^\perp\bn_k\wedge\bH_k\big)\:.
\ees
We suppose that
\be\label{d1}
div\,Q_k\:\longrightarrow\:0\qquad\text{strongly in}\quad \big(W^{2,2}\cap W^{1,\infty}\big)'(\di)\:.
\ee
and
\be\label{d2}
\big\Vert\bp_k\Vert_{W^{2,2}\cap W^{1,\infty}}\;\le\;C\qquad\text{uniformly, for some constant $C>0$\;.}
\ee

\medskip
\noindent
We begin our study with an elementary result. 

\begin{Lma}
\label{conv1}
There holds
\bes
\bp_k\;div\;Q_k\:\longrightarrow\:0\qquad\text{in}\qquad\mathcal{D}'(\di)\:.
\ees
\end{Lma}
$\textbf{Proof.}$ For notational convenience, we set
\bes
X\,:=\,W^{2,2}\,\cap\,W^{1,\infty}\:.
\ees\\[-3ex]
Let $u$ be an arbitrary element of $\,\mathcal{D}=C_0^\infty(\di)$\,. The Sobolev embedding theorem guarantees that the elements $\bp_k$ of $X$ are H\"older continuous on $\overline{\di}$. Since in addition $X$ is an intersection of Sobolev spaces, it is clear that $\,\bp_k u\,$ is an element of $X_0$ (i.e. an element of $X$ with null trace on the boundary of the unit disk). 
More precisely, from (\ref{d2}), 
\be\label{dd1}
\big\Vert\bp_k\,u\big\Vert_{X_0}\;\lesssim\;\big\Vert\bp_k\big\Vert_{X}\,\big\Vert u\big\Vert_{\mathcal{D}}\;\le\;C\,\big\Vert u\big\Vert_{\mathcal{D}}\:.
\ee
Setting $\,T_k:=div\;Q_k\,$, we can make sense of $\bp_k\,T_k$ as a distribution via
\bes
\big\langle \bp_k\,T_k\,,u\big\rangle_{\mathcal{D}',\mathcal{D}}\;:=\;\big\langle T_k\,,\bp_k\,u\big\rangle_{X_0',X_0}\:.
\ees
It follows immediately from (\ref{d1}) and (\ref{dd1}) that $\bp_k\,T_k$ converges to zero in the sense of distributions.\\[-3ex]

$\hfill\blacksquare$\\

\noindent
By hypothesis, $\,T_k:=div\,Q_k\,$ is a bounded linear functional on $X$ (with $X$ as in the proof of the previous lemma). Using (\ref{c1}), we see that $T_k$ belongs to $W^{-2,(2,\infty)}$, which is the dual space of $\,W_0^{2,(2,1)}$. Accordingly, Proposition~\ref{charsob} grants the existence of an element $P_k$ of $\,\mathbb{R}^2\otimes W^{-1,(2,\infty)}\,$ satisfying
\be\label{d3}
T_k\;=\;div\,P_k
\ee
and
\be\label{d4}
\Vert P_k \Vert_{\mathbb{R}^2\otimes W^{-1,(2,\infty)}}\;\le\;\Vert T_k\Vert_{W^{-2,(2,\infty)}}\,+\,\delta_k\:,
\ee

\medskip
\noindent
for some positive constant $\delta_k$ arbitrarily chosen. \\[1ex]
We next establish a useful result. 
\begin{Lma}
\label{conv2}
There holds
\bes
P_k\,\bp_k\:\longrightarrow\:0\qquad\text{in}\qquad\mathbb{R}^2\otimes\mathcal{D}'(\di)\:.
\ees
\end{Lma}
$\textbf{Proof.}$ From (\ref{d4}) and the inclusion
\bes
X\;:=\;W^{2,2}\cap W^{1,\infty}\;\supset\;W_0^{2,(2,1)}\:,
\ees
there holds
\bes
\Vert P_k\Vert_{\mathbb{R}^2\otimes W^{-1,(2,\infty)}}\:\le\:\Vert T_k\Vert_{X'}\,+\,\delta_k\:.
\ees

\noindent
Owing to (\ref{d1}), the sequence $(P_k)$ converges strongly to zero in $\,\mathbb{R}^2\otimes W^{-1,(2,\infty)}$. 
Moreover, from (\ref{d2}) and (\ref{c11}), it follows that the sequence $(\bp_k)$ is uniformly bounded in $\,W^{1,(2,1)}$. The desired statement may now be reached by repeating mutatis mutandis the proof of Lemma~\ref{conv1}, and letting $\delta_k$ tend to zero.\\[-4ex]

$\hfill\blacksquare$

\medskip
\noindent
The identity (\ref{d3}) yields
\be\label{d5}
div\,\big(Q_k-P_k\big)\;=\;0\:.
\ee
By definition, the difference $(Q_k-P_k)$ lies in $\,\mathbb{R}^2\otimes W^{-1,(2,\infty)}$. As proved in Lemma~\ref{A-1}, the equation (\ref{d5}) implies the existence of an element $\bL_k$ in $L^{2,\infty}$ satisfying
\be\label{d6}
Q_k\,-\,P_k\;=\;\nabla^{\perp}\bL_k\:.
\ee
Since $(\bp_k)$ is uniformly bounded in $W^{2,2}\cap\,W^{1,\infty}$, it follows that the sequences $(\bH_k)$, $(\bn_k)$, and thus $(Q_k)$ are uniformly bounded, respectively in $L^{2}$, $W^{1,2}$, and $\mathbb{R}^2\otimes W^{-1,(2,\infty)}$. Likewise, we saw that $(P_k)$ is uniformly bounded $\mathbb{R}^2\otimes W^{-1,(2,\infty)}$. Hence, the sequence $\,(\nabla^{\perp}\bL_k)$\, is uniformly bounded in $\,\mathbb{R}^2\otimes W^{-1,(2,\infty)}$.

\begin{Lma}
\label{conv3}
There holds
\bes
\left\{\begin{array}{rcl}
\nabla\bp_k\cdot\nabla^\perp\bL_k&\longrightarrow&0\\[1.25ex]
\nabla\bp_k\wedge\big(\nabla^\perp \bL_k\,+\,2\,\nabla\bH_k\big)&\longrightarrow&0
\end{array}\right.\qquad\text{in}\qquad\mathcal{D'}(\di)\:.
\ees
\end{Lma}
$\textbf{Proof.}$ Observe that
\bes
P_k\cdot\nabla\bp_k\;=\;div\,\big(P_k\,\bp_k\big)\,-\,\bp_k\;div\,P_k\;=\;div\,\big(P_k\,\bp_k\big)\,-\,\bp_k\;div\,Q_k\:.
\ees
Whence, the results of Lemma~\ref{conv1} and Lemma~\ref{conv2} show that
\be\label{y0y}
\nabla\bp_k\cdot P_k\:\longrightarrow\:0\qquad\text{in}\qquad\mathcal{D}'(\di)\:.
\ee
From (\ref{b1}) in the Appendix, we know that 
\bes
\nabla\bp_k\cdot Q_k\;=\;0\qquad\text{and}\qquad\nabla\bp_k\wedge Q_k\;=\;-\,2\,\nabla\bp_k\wedge\nabla\bH_k\:.
\ees
The desired statement now ensues by combining (\ref{d6}) and (\ref{y0y}).\\[-3ex]

$\hfill\blacksquare$\\

\noindent
From the characterization provided in Proposition~\ref{charsob}, we can always arrange for the $L^{2,\infty}$-norm of $\bL_k$ to be as close as we please to the $W^{-1,(2,\infty)}$-norm of $\nabla^{\perp}\bL_k$. In particular, the sequence $(\bL_k)$ is uniformly bounded in $L^{2,\infty}$. 
We may thus extract a weak* convergent subsequence with
\be\label{d61}
\left\{\begin{array}{rclcccl}
\bL_{k'}&\stackrel{*}{-\!\!\!\!-\!\!\!\rightharpoondown}&\bL&\qquad&\text{in}&&L^{2,\infty}\\[1ex]
\nabla^\perp\bL_{k'}&\stackrel{*}{-\!\!\!\!-\!\!\!\rightharpoondown}&g&\quad&\text{in}&& \mathbb{R}^2\otimes W^{-1,(2,\infty)}\:,
\end{array}\right.
\ee
with $\,g=\nabla^{\perp}\bL$\, in the sense of distributions. \\[1ex]
As the sequence $(\bp_k)$ is uniformly bounded in $\,W^{2,2}\cap W^{1,\infty}$\, by hypothesis, the Banach-Alaoglu theorem implies that there exists a subsequence $(\bp_{k'})$ converging weak* in $\,W^{2,2}\cap W^{1,\infty}$\, to some element $\bp$. In turn, the compact embeddings provided by the Rellich-Kondrachov theorem (\ref{c11}) enables us to further extract a subsequence, still denoted $(\bp_{k'})$, satisfying the strong convergences
\be\label{d12}
\left\{\begin{array}{rclcccl}
\bp_{k'}&\longrightarrow&\bp&\qquad&\text{in}&&\bigcap_{p<\infty}W^{1,p}\\[1ex]
\nabla\bp_{k'}&\longrightarrow&\nabla\bp&\quad&\text{in}&&\mathbb{R}^2\otimes\bigcap_{p<\infty}L^{p}\:.
\end{array}\right.
\ee
The convergences (\ref{d61}) and (\ref{d12}) yield

\begin{Lma}
\label{conv4}
\bes
\left\{\begin{array}{lclc}
\nabla\bp_{k'}\cdot\nabla^\perp\bL_{k'}&\longrightarrow&\nabla\bp\,\cdot\nabla^\perp\bL\\[1.5ex]
\nabla\bp_{k'}\wedge\nabla^\perp\bL_{k'}&\longrightarrow&\nabla\bp\,\wedge\nabla^\perp\bL
\end{array}\right.\qquad\text{in}\qquad\mathcal{D}'(\di)\:.
\ees
\end{Lma}
$\textbf{Proof.}$ We shall only establish the first convergence, the second one being obtained {\it mutatis mutandis}. Owing to the general identity
\bes
div\big(a\,\nabla^{\perp}b)\;=\;\nabla a\cdot\nabla^\perp b\:,
\ees
it suffices to show that
\bes
\bp_{k'}\,\nabla^\perp\bL_{k'}\:\longrightarrow\:\bp\;\,\nabla^\perp\bL\qquad\text{in}\qquad\mathcal{D}'(\di)\:.
\ees
This is what we shall do. For convenience, we set $\,Y:=\mathbb{R}^2\otimes W^{1,(2,1)}$. Let $g$ be as in (\ref{d61}), and $u$ be an arbitrary test-function in $\mathbb{R}^2\otimes\mathcal{D}$. Clearly, $\bp\,u$ and $\bp_{k'}u$ are elements of $Y_0$ (i.e. elements of $Y$ with null trace on the boundary of the unit-disk), with an estimate analogous to (\ref{dd1}). Note that
\begin{eqnarray*}
\Big\langle \bp_{k'}\nabla^\perp\bL_{k'}\,-\,\bp\,g\,,u\Big\rangle_{\mathcal{D}',\mathcal{D}}&&\\[1ex]
&&\hspace{-2.5cm}\equiv\:\;\Big\langle\nabla^\perp\bL_{k'}\,,\,u\,\big(\bp_{k'}-\bp\big)\Big\rangle_{Y_0',Y_0}\;+\;\Big\langle\nabla^\perp\bL_{k'}-g\,,\,u\,\bp\Big\rangle_{Y_0',Y_0}\:.
\end{eqnarray*}
Whence, from the convergences (\ref{d61}) and (\ref{d12}), the desired result follows. \\[-3ex]

$\hfill\blacksquare$

\medskip

As explained in Theorem 3.3.8 from \cite{Hel}, the fact that\footnote{where $K_k$ denotes the Gaussian curvature associated with the immersion $\bp_k$.}
\bes
\Delta\la_k\;=\;-\,\text{e}^{2\la_k}K_k
\ees
is an element of the Hardy space $\mathcal{H}^1$ implies that $(\la_k)$ is a sequence of elements in $W^{1,(2,1)}$, uniformly bounded in norm by a constant depending only upon the uniform bound on $\,\big(\Vert\bp_k\Vert_{W^{2,2}\cap W^{1,\infty}}\big)$. Owing to the Rellich-Kondrachov theorem, we may extract a subsequence $(\la_{k'})$ satisfying
\be\label{d13}
\la_{k'}\:\longrightarrow\:\la\qquad\text{in}\qquad\bigcap_{p<\infty}L^p\:,
\ee
for some $\la$ in the suitable space. Recall that for $j\in\{1,2\}$, there holds
\bes
\text{e}^{2\la_k}\;=\;\big\vert\partial_{x^j}\bp_k\big\vert^2\qquad\text{and thus}\qquad\partial_{x^i}\big(\text{e}^{2\la_k}\big)\;=\;2\,\partial_{x^j}\bp_k\cdot\partial_{x^jx^i}\bp_k\:.
\ees
Accordingly, the sequence $(\text{e}^{2\la_k})$ is uniformly bounded in $W^{1,2}\cap L^{\infty}$. From this, and the boundedness of $\la_k$, the same is true about the sequences $(\text{e}^{\pm\la_k})$. Combining this to the convergences (\ref{d12}) and (\ref{d13}), we find the following strong convergences
\be\label{d14}
\left\{\begin{array}{rcl}
\text{e}^{\pm\la_{k'}}&\longrightarrow&\text{e}^{\pm\la_{k}}\\[1ex]
(\bej)_{k'}\;:=\;\text{e}^{-\la_{k'}}\partial_{x^j}\bp_{k'}&\longrightarrow&\text{e}^{-\la}\partial_{x^j}\bp\;=:\;\bej\end{array}\right.\quad\text{in}\quad\bigcap_{p<\infty}L^p\:.
\ee
By definition, from the uniform boundedness of the sequence $(\bp_k)$ in $W^{2,2}\cap W^{1,\infty}$ and of that of the sequence $(\bn_k)$ in $W^{1,2}$, we see that the sequence $(\bH_k)$ is uniformly bounded in $L^2$. Accordingly, we can extract a weakly convergent subsequence:
\be\label{d15}
\bH_{k'}\:-\!\!\!\!-\!\!\!\rightharpoondown\:\bH\qquad\text{in}\qquad L^2\:.
\ee
Altogether, (\ref{d14}) and (\ref{d15}) show that
\be\label{d16}
(\bex)_{k'}\wedge(\bey)_{k'}\wedge\bH_{k'}\:-\!\!\!\!-\!\!\!\rightharpoondown\:\bex\wedge\bey\wedge\bH\qquad\text{in}\qquad \bigcap_{1\le p<2}L^p\:.
\ee
Equation (II.42) established in \cite{Riv2} states
\bes
\star\,(\bn_{k'}\res\bH_{k'})\;=\;(-1)^{m-1}\,(\bex)_{k'}\wedge(\bey)_{k'}\wedge\bH_{k'}\:.
\ees
Consequently, (\ref{d16}) is tantamount to
\bes
\star\,(\bn_{k'}\res\bH_{k'})\:-\!\!\!\!-\!\!\!\rightharpoondown\:\star\,(\bn\res\bH)\qquad\text{in}\qquad \bigcap_{1\le p<2}L^p\:.
\ees
Combining this to (\ref{d12}) shows that
\be\label{d17}
\star\,(\bn_{k'}\res\bH_{k'})\,\res\,\nabla^\perp\bp_{k'}\:-\!\!\!\!-\!\!\!\rightharpoondown\:\star\,(\bn\,\res\bH)\,\res\,\nabla^\perp\bp\qquad\text{in}\qquad\bigcap_{1\le p<2}L^p\:.
\ee
It takes little effort to verify that
\bes
div\,(u\,\res \nabla^{\perp}v)\;=\;\nabla u\,\res\,\nabla^{\perp}v\:.
\ees
holds in general for suitable $u$, $v$. In particular, (\ref{d17}) yields
\bes
\nabla\big(\star(\bn_{k'}\res\bH_{k'})\big)\,\res\,\nabla^\perp\bp_{k'}\:\longrightarrow\:\nabla\big(\star(\bn\,\res\bH)\big)\,\res\,\nabla^\perp\bp\qquad\text{in}\qquad\mathcal{D}'(\di)\:.
\ees
Using identity (II.49) from \cite{Riv2}, the latter is equivalent to
\be\label{d18}
\nabla\bp_{k'}\wedge\nabla\bH_{k'}\:\longrightarrow\:\nabla\bp\,\wedge\nabla\bH\qquad\text{in}\qquad\mathcal{D}'(\di)\:.
\ee

Finally, bringing altogether Lemma~\ref{conv3}, Lemma~\ref{conv4}, and (\ref{d18}), we conclude as desired that, in the sense of distributions,
\bes
\left\{\begin{array}{lcl}
\nabla\bp\,\cdot\nabla^\perp\bL&=&0\\[1ex]
\nabla\bp\wedge\nabla^\perp\bL&=&-\,2\,\nabla\bP\wedge\nabla\bH\:.
\end{array}\right.
\ees

\subsection{Proof of Theorem~\ref{th-cw}}

As done in the work \cite{Riv2}, we introduce the functions $S$ and $\bR$ such that
\be\label{r2}
\left\{\begin{array}{lcl}
\nabla S&:=&\nabla\bp\cdot\bL\\[1ex]
\nabla\bR&:=&\nabla\bp\wedge\bL\,+\,2\,\nabla^\perp\bp\wedge\bH\:.
\end{array}\right.
\ee
Clearly, $S$ and $\bR$ are defined on $\di$ up to an unimportant additive constant.\\
By hypothesis, $\bL$ is an element of $L^{2,\infty}(\di)$, while $\bp$ is Lipschitz. Accordingly, $\nabla S$ and $\nabla\bR$ belong to $L^{2,\infty}(\di)$ and to $\R^2\otimes L^{2,\infty}(\di)$ respectively. This observation, and the particular structure of the right-hand side of the system (\ref{r2}) will enable us to deduce the regularity result we are seeking to obtain.

\subsubsection{Regularity}

It is shown in \cite{Riv2} that $S$ and $\bR$ satisfy the system
\be\label{r3}
\left\{\begin{array}{lcl}
\Delta S&=&\big(\nabla\star\bn\big)\cdot\nabla^\perp\bR\\[1ex]
\Delta\bR&=&(-1)^n\star\big(\nabla\bn\bul\nabla^\perp\bR\big)\,-\,(\nabla\star\bn)\,\nabla^\perp S\:.
\end{array}\right.
\ee
The advantage of these equations lies essentially in their right-hand sides comprising only Jacobians. This peculiar feature will enable us to apply the techniques of integration by compensation, more precisely Wente-type estimates. We recall two results which shall be of use to our proof. They are due in parts to contributions by Wente \cite{Wen}, Tartar \cite{Ta2}, Coifman {\it et al.} \cite{CLMS}, Bethuel \cite{Be}, and H\'elein \cite{Hel}.  
\begin{Lma}
\label{reg1}
Let $\Om$ be an open subset of $\mathbb{R}^2$ with $C^2$-boundary. Suppose that $a$ and $b$ are elements of $W^{1,2}(\Om)$ and of $W^{1,(2,\infty)}(\Om)$, respectively. If $u$ satisfies
\bes
\left\{\begin{array}{rclcl}
\Delta u&=&\nabla a\cdot\nabla^\perp b&\quad&\text{in}\:\:\:\Om\\[1ex]
u&=&0&\quad&\text{on}\:\:\:\partial\Om\:,
\end{array}\right.
\ees
then $\nabla u$ belongs to the space $L^{2}(\Om,\mathbb{R}^2)$ with the estimate
\bes
\big\Vert\nabla u\big\Vert_{L^{2}(\Om,\mathbb{R}^2)}\:\:\lesssim\:\:\big\Vert\nabla a\big\Vert_{L^{2}(\Om,\mathbb{R}^2)}\,\big\Vert\nabla b\big\Vert_{L^{2,\infty}(\Om,\mathbb{R}^2)}\:,
\ees
up to a multiplicative constant depending only on $\Om$.
\end{Lma}
The interested reader will find the proof of this result and further variations on the same theme in \cite{Hel} (Theorem 3.4.5).\\
Another result (Theorem 3.4.1 in \cite{Hel}) which will be useful to us is the following.

\begin{Lma}
\label{reg2}
Let $\Om$ be an open subset of $\mathbb{R}^2$ with $C^1$-boundary. Suppose that $a$ and $b$ are elements of $W^{1,2}(\Om)$. If $u$ satisfies
\bes
\left\{\begin{array}{rclcl}
\Delta u&=&\nabla a\cdot\nabla^\perp b&\quad&\text{in}\:\:\:\Om\\[1ex]
u&=&0&\quad&\text{on}\:\:\:\partial\Om\:,
\end{array}\right.
\ees
then $u$ belongs to the space $\,W^{1,(2,1)}(\Om)\subset C^0(\Om)$ with the estimate
\bes
\big\Vert \nabla u\big\Vert_{L^{2,1}(\Om)}\:\:\lesssim\:\:\big\Vert\nabla a\big\Vert_{L^{2}(\Om,\mathbb{R}^2)}\,\big\Vert\nabla b\big\Vert_{L^{2}(\Om,\mathbb{R}^2)}\:,
\ees
up to a multiplicative constant depending only on $\Om$.
\end{Lma}

Geared with these results, we are prepared to start our proof. Let us define
\be\label{r6}
S\;=\;S_0\,+\,S_1\hspace{1cm}\text{and}\hspace{1cm}\bR\;=\;\bR_0\,+\,\bR_1\:,
\ee
where the new variables, in accordance with (\ref{r3}), satisfy
\be\label{r4}
\left\{\begin{array}{rclcrclcl}
\Delta S_0&=&0&,\qquad&\Delta\bR_0&=&\vec{0}&\quad&\text{in}\:\:\:\di\\[1ex]
S_0&=&S&,\qquad&\bR_0&=&\bR&\quad&\text{on}\:\:\:\partial\di\:,
\end{array}\right.
\ee
and
\be\label{r5}
\left\{\begin{array}{rclcrclcl}
\Delta S_1&=&\Delta S&,\qquad&\Delta\bR_1&=&\Delta\bR&\quad&\text{in}\:\:\:\di\\[1ex]
S_1&=&0&,\qquad&\bR_1&=&0&\quad&\text{on}\:\:\:\partial\di\:.
\end{array}\right.
\ee\\
We saw above that $S$ and $\bR$ are elements of $W^{1,(2,\infty)}(\di)$, while $\star\,\bn$ belongs to $W^{1,2}(\di)$. Since the right-hand sides of (\ref{r3}), and thus of (\ref{r5}), comprise only Jacobians, Lemma~\ref{reg1} may be called upon so as to produce the estimates
\be\label{r7}
\left\{\begin{array}{lcl}
\big\Vert\nabla S_1\big\Vert_{L^{2}(\di)}&\lesssim&\eps\,\big\Vert\nabla\bR\big\Vert_{L^{2,\infty}(\di)}\\[2ex]
\big\Vert\nabla \bR_1\big\Vert_{L^{2}(\di)}&\lesssim&\eps\,\Big(\big\Vert\nabla S\big\Vert_{L^{2,\infty}(\di)}\,+\,\big\Vert\nabla\bR\big\Vert_{L^{2,\infty}(\di)}\Big)\:,
\end{array}\right.
\ee
up to unimportant multiplicative constants depending only on $\di$. Here $\eps$ denotes the (adjustable) upper bound on the energy:
\be\label{nrj}
\big\Vert\nabla\bn\big\Vert_{L^2(\di)}\:\:\le\:\:\eps\:.
\ee

\noindent
Let us fix once and for all some point $p\in D_{1/2}(0)$, some $\,0<k<1\,$, and some radius $\,0<r<1/2$. The flat disk $D_{kr}(p)$ of radius $kr$ and centered on the point $p$ is properly contained in the unit-disk $\di$. Since $S_0$ and $\bR_0$ are harmonic and satisfy (\ref{r4}), it is clear that $\nabla S_0$ and $\nabla\bR_0$ are square-integrable on $D_{kr}(p)$. From (\ref{r7}), it follows that $\nabla S$ and $\nabla\bR$ are likewise square-integrable on $D_{kr}(p)$. Because $\star\,\bn$ is also an element of $W^{1,2}(D_{kr}(p))$, we are now in position to apply Lemma~\ref{reg2} to the system (\ref{r5}), thereby obtaining
\bes
\left\{\begin{array}{lcl}
\big\Vert\nabla S_1\big\Vert_{L^{2,1}(D_{kr}(p))}&\lesssim&\eps\,\big\Vert\nabla\bR\big\Vert_{L^{2}(D_{kr}(p))}\\[2ex]
\big\Vert\nabla \bR_1\big\Vert_{L^{2,1}(D_{kr}(p))}&\lesssim&\eps\,\Big(\big\Vert\nabla S\big\Vert_{L^{2}(D_{kr}(p))}\,+\,\big\Vert\nabla\bR\big\Vert_{L^{2}(D_{kr}(p))}\Big)\:,
\end{array}\right.
\ees
up to some unimportant multiplicative constants involving $\eps$ (they are independent of $r$ by homogeneity).\\[.5ex]
In particular, since $L^{2,1}\subset L^2$ and $k<1$, the latter gives that for some positive constant $C_0$, there holds
\be\label{stuff1}
E_{D_{kr}(p)}(S_1,\bR_1)\:\le\:C_0\,\eps\,E_{D_{r}(p)}(S,\bR)\:,
\ee
where, for notational convenience, we have set
\bes
E_{D_\rho(p)}(u,\bv)\::=\:\big\Vert u\big\Vert_{L^2(D_\rho(p))}\,+\,\big\Vert\bv\big\Vert_{L^2(D_\rho(p))}\:.
\ees
On the other hand, since $S_0$ and $\bR_0$ are harmonic, a classical ``monotonicity" result (cf. Theorem 3.3.12 in \cite{Hel}, and \cite{Gia}) applied to (\ref{r4}) yields that
\be\label{stuff2}
E_{D_{kr}(p)}(S_0,\bR_0)\:\le\:k\,E_{D_{r}(p)}(S_0,\bR_0)\:\le\:k\,E_{D_{r}(p)}(S,\bR)\:.
\ee
Via combining altogether (\ref{stuff1}) and (\ref{stuff2}), we obtain
\be\label{stuff3}
E_{D_{kr}(p)}(S,\bR)\:\le\:(C_0\,\eps+k)\,E_{D_{r}(p)}(S,\bR)\:.
\ee
We have the freedom to adjust the positive parameter $\eps$ as we please. Because $k\in(0,1)$, we may in particular arrange for the constant $(C_0\,\eps+k)$ to be smaller than 1. Then, iterating (\ref{stuff3}), we infer the existence of some $\gamma\in(0,1)$ such that
\be\label{stuff4}
E_{D_\rho(p)}\:\le\:C\,\rho^\gamma
\ee
holds for all $\,0<\rho<1/2$, all points $p\in D_{1/2}(0)$, and some constant $C>0$. With the help of the Poincar\'e inequality, the estimate (\ref{stuff4}) can be used to show that $S$ and $\bR$ are locally H\"older continuous (see \cite{Gia}). We are however interested in another corollary of (\ref{stuff4}). We denote the maximal function by
\bes
M_{2-\gamma}F(x)\::=\:\sup_{\rho>0}\;\rho^{-\gamma}\int_{D_\rho(x)}|F(y)|\,dy\:.
\ees
Going back to (\ref{r3}), owing to (\ref{stuff4}), it follows that 
\bes
M_{2-\gamma}\chi_{D_{1/2}(0)}\Delta S(p)\:\le\:\eps\,\sup_{0<\rho<\frac{1}{2}}\;\rho^{-\gamma}\,E_{D_\rho(p)}\:\le\:C\,\eps\:,\qquad\forall\:\:p\in D_{1/2}(0)\:.
\ees
The same estimate clearly holds with $\bR$ in place of $S$. Accordingly,
\bes
\left\{\begin{array}{rcl}
\big\Vert M_{2-\gamma}\chi_{D_{1/2}(0)}\Delta S\big\Vert_{L^{\infty}(D_{1/2}(0))}&<&\infty\\[2ex]
\big\Vert M_{2-\gamma}\chi_{D_{1/2}(0)}\Delta\bR\big\Vert_{L^{\infty}(D_{1/2}(0))}&<&\infty\:.
\end{array}\right.
\ees
Moreover, it is clear that $\Delta S$ and $\Delta\bR$ belong to $L^1(D_{1/2}(0))$, since $S$, $\bR$, and $\star\,\bn$ are elements of $W^{1,2}(D_{1/2}(0))$. We may thus call upon Proposition 3.2 in \cite{Ad2} to deduce that 
\bes
\dfrac{1}{|x|}*\chi_{D_{1/2}(0)}\Delta S\qquad\text{and}\qquad \dfrac{1}{|x|}*\chi_{D_{1/2}(0)}\Delta\bR
\ees
belong to $\,L^{q,\infty}(D_{1/2}(0))$, with $\,q=\dfrac{2-\gamma}{1-\gamma}\,$.\\[1ex]
A classical estimate about Riesz kernels states that in general, there holds
\bes
|\nabla u|(p)\:\le\:C_1\,\dfrac{1}{|x|}*\chi_{D_{1/2}(0)}\Delta u\,+\,C_2\:,\qquad\forall\:\:p\in D_{1/4}(0)\:,
\ees
for two constant $C_1$ and $C_2$. Hence, we infer that $\nabla S$ and $\nabla\bR$ are elements of $L^{q,\infty}(D_{1/4}(0))$, with $q$ as above. In particular, because $\gamma\in(0,1)$, it follows that
\be\label{xx0}
\nabla S\,\in\,L^{2+\delta}(D_{1/4}(0),\R)\qquad\text{and}\qquad \nabla\bR\,\in\,L^{2+\delta}(D_{1/4}(0),\R^2)\:,
\ee
or some $\delta>0$.\\[-1ex]

Now that we have this result at our disposal, we are ready to implement our regularity proof. It involves a bootstrapping argument, which will be performed on the following identity, established in Lemma~\ref{A-4}\,:
\be\label{r1}
-\,\Delta\bp\;=\;\nabla\bR\,\bul\nabla^\perp\bp\;+\,\nabla S\;\nabla^\perp\bp\:.
\ee
By hypothesis, $\bP$ is a Lipschitz function. Combining (\ref{xx0}) and (\ref{r1}) thus shows that 
\be\label{xx1}
\nabla\bP\,\in\,W^{1,2+\delta}(D_{1/4}(0),\R^2)\:.
\ee

\noindent
Recall next that
\be\label{r8}
\star\,\bn\;=\;\text{e}^{-2\la}\,\partial_{x^1}\bp\,\wedge\partial_{x^2}\bp\qquad\text{with}\qquad 2\,\text{e}^{2\la}\;=\;\big|\nabla\bp\big|^2\:.
\ee
It is known that $\la$ belongs to $W^{1,(2,1)}(\di)$. Indeed, in \cite{MS} the authors show that
\bes
\Delta\la\;=\;-\,\text{e}^{2\la}\,K\:,\qquad\text{where $K$ is the Gaussian curvature}\:,
\ees 
is an element of the Hardy space $\mathcal{H}^1$, thereby implying that $\nabla\la$ lies in  the Lorentz space $L^{2,1}$ (cf. Theorem 3.3.8 in \cite{Hel}). Accordingly, $\text{e}^{\pm2\la}$ are bounded from above and below.\\
Suppose now that $\nabla\bP$ belongs to $W^{1,a}$, for some $a>2$. Since by hypothesis $\bP$ is Lipschitz, the second equation in (\ref{r8}) implies that $e^{2\la}$ is an element of $W^{1,a}\cap L^{\infty}$. In turn, because $e^{\pm\la}$ are bounded from above and below, the first equation in (\ref{r8}) yields that $\nabla\star\bn$ belongs to $L^{a}$. In particular, (\ref{xx1}) shows that
\be\label{xx2}
\nabla\star\bn\,\in\,L^{2+\delta}(D_{1/4}(0),\R^2)\:.
\ee

\medskip
\noindent
For the reader's convenience, we recall
\bes
\left\{\begin{array}{lcl}
\Delta S&=&\big(\nabla\star\bn\big)\cdot\nabla^\perp\bR\\[1ex]
\Delta\bR&=&(-1)^n\star\big(\nabla\bn\bul\nabla^\perp\bR\big)\,-\,(\nabla\star\bn)\,\nabla^\perp S\:.
\end{array}\right.
\ees
Accounting for (\ref{xx0}) and (\ref{xx2}) in this system gives
\bes
\nabla S\,\in\,W^{1,1+\delta/2}\,\subsetneq\,L^{2+\delta}\:.
\ees
Naturally, the same statement holds with $\bR$ in place of $S$. Comparing the latter to (\ref{xx0}) reveals that the regularity has been improved. The process may thus be repeated a finite of number of times, until we reach that $\nabla S$ and $\nabla\bR$ are continuous. The exact same procedure as above will then yield that $\bP$ belongs to $C^1$, and eventually that it is smooth, thereby concluding the proof.

\subsubsection{Conformal Willmore Equation}

We open our derivations by introducing some notation. Let $z=x^1+ix^2$ and $z^*$ be its complex conjugate. We set
\bes
\bez\;:=\;\text{e}^{-\la}\,\partial_z\bp\;=\;\dfrac{1}{2}\,\big(\bex\,-\,i\,\bey\big)\hspace{.75cm}\text{and}\hspace{.75cm}\bezz\;:=\;\text{e}^{-\la}\,\partial_{z^*}\bp\;=\;\dfrac{1}{2}\,\big(\bex\,+\,i\,\bey\big)\:.
\ees
Since $\bp$ is conformal, there holds
\bes
\pj\bp\cdot\pk\bp\;=\;\text{e}^{2\la}\,\delta_{jk}\hspace{1cm}\text{for}\qquad(j,k)\,\in\,\{1,2\}\:,
\ees
and thus
\bes
\pa\bp\cdot\pb\bp\;=\;\dfrac{1}{2}\,\text{e}^{2\la}\,\delta_{ab^*}\hspace{1cm}\text{for}\qquad(a,b)\,\in\,\{z,z^{*}\}\:.
\ees
From this, it follows easily that for any triple $\,(a,b,c)\in\{z,z^*\}\,$ there holds
\begin{eqnarray*}
\pa\bp\cdot\partial_{bc}\bp&\equiv&\delta_{ca}\big(\pa\bp\cdot\partial_{ab}\bp\big)\,+\,\delta_{cb}\big(\pa\bp\cdot\partial_{bb}\bp\big)\\[1ex]
&=&\dfrac{1}{2}\,\big(\delta_{ca}-\delta_{cb}\big)\pb\big|\pa\bp\big|^2\;+\;\delta_{cb}\,\pb\big(\pa\bp\cdot\pb\bp\big)\\[1ex]
&=&\delta_{cb}\,\delta_{ab^*}\,\text{e}^{2\la}\,\partial_b\la\:.
\end{eqnarray*}
Hence,
\begin{eqnarray}\label{a1}
\bbe_a\cdot\pb\bbe_c&\equiv&\text{e}^{-2\la}\,\big(\pa\bp\cdot\partial_{bc}\bp\,-\,(\pb\la)\,\pa\bp\cdot\partial_c\bp\big)\nonumber\\[1ex]
&=&\dfrac{1}{2}\,\big(2\,\delta_{cb}\,\delta_{ab^*}\,-\,\delta_{ac^*}\big)\,\pb\la\:.
\end{eqnarray}
Observe furthermore that
\be\label{a2}
\bbe_a\cdot\bbe_b\;=\;\dfrac{1}{2}\,\delta_{ab^*}\:.
\ee
Thus, combining (\ref{a1}) and (\ref{a2}) gives
\be\label{a3}
\pb\bbe_c\;=\;\partial_b\la\sum_{a\in\{z,z^*\}}\big(2\delta_{cb}\delta_{ab}-\delta_{ac}\big)\,\bbe_a\;+\;\pro\big(\pb\bbe_c\big)\:.
\ee
Next, we have
\begin{eqnarray*}
\bn_\al\cdot\pzz\bez&=&-\,\bez\cdot\pzz\bn_\al\:\:=\:\:-\,\dfrac{1}{4}\big(\bex-i\,\bey)(\px+i\,\py)\,\bn_\al\\[1ex]
&=&-\,\dfrac{\text{e}^\la}{4}\,\big(h^\al_{11}+h^\al_{22}\big)\:\:=\:\:-\,\dfrac{\text{e}^\la}{2}\,H^\al\:,
\end{eqnarray*}
and similarly
\bes
\bn_\al\cdot\pzz\bezz\;=\;\dfrac{\text{e}^\la}{2}\,H^\al_0\:.
\ees

\medskip
\noindent
Accordingly, we may now deduce from (\ref{a3}) that
\be\label{a4}
\left\{\begin{array}{lcl}
\pzz\bez&=&-\,\big(\pzz\la\big)\,\bez\;+\;\dfrac{\text{e}^\la}{2}\,\vec{H}\\[2.5ex]
\pzz\bezz&=&\big(\pzz\la\big)\,\bezz\;+\;\dfrac{\text{e}^\la}{2}\,\vec{H}_0\:.
\end{array}\right.
\ee
Note that we can also easily obtain from the above computations the identity
\be\label{a5}
\pzz\bn_\al\:=\:-\,{\text{e}^\la}\,\big(H_0^\al\,\bez\,+\,H^\al\,\bezz\big)\,+\,\pro\big(\pzz\bn_\al\big)\:.
\ee
These expressions shall come helpful in the sequel. \\

Suppose now that $\bp$ is smooth and satisfies for some $\bL$ the system
\be\label{a11}
\left\{\begin{array}{lcl}
\nabla\bp\cdot\nablap\bL&=&0\\[1ex]
\nabla\bp\wedge\big(2\,\nabla\bH+\nablap\bL\big)&=&{0}\:.
\end{array}\right.
\ee
Since $\,\bbe_j=\text{e}^{-\la}\,\pj\bp\,$, the first equation in the system is equivalent to
\bes
\bex\cdot\py\bL\;=\;\bey\cdot\px\bL\:.
\ees
Whence, we deduce that $\bL$ satisfies
\be\label{a12}
\nabla^\perp\bL\;=\;\left(\begin{array}{cc}a&b\\[1ex]c&-\,a\end{array}\right)\left(\begin{array}{c}\bex\\[1ex]\bey\end{array}\right)\,+\,\left(\begin{array}{c}p^\al\\[1ex]q^\al\end{array}\right)\bn_\al\:,
\ee
for some suitable coefficients $a$, $b$, $c$, $p^\al$, and $q^\al$. \\[1ex]
Substituting this form in the second equation from (\ref{a11}) gives
\begin{eqnarray*}
&&\bex\wedge\Big[\,2\,\big(\bex\cdot\px\bH\big)\,\bex\,+\,2\,\big(\bey\cdot\px\bH\big)\,\bey\,+\,2\,\big(\bn_\al\cdot\px\bH\big)\,\bn_\al\\
&&\hspace{6cm}+\;a\,\bex\,+\,b\,\bey\,+\,p^\al\,\bn_\al\Big]\\[1ex]
&=&-\,\bey\wedge\Big[\,2\,\big(\bex\cdot\py\bH\big)\,\bex\,+\,2\,\big(\bey\cdot\py\bH\big)\,\bey\,+\,2\,\big(\bn_\al\cdot\py\bH\big)\,\bn_\al\\
&&\hspace{6cm}+\;c\,\bex\,-\,a\,\bey\,+\,q^\al\,\bn_\al\Big]\:,
\end{eqnarray*}
thereby yielding
\bes
c-b\;=\;2\,(\bey\cdot\px-\bex\cdot\py)\bH\:,\qquad p^\al\;=\;-\,2\,\bn_\al\cdot\px\bH\:,\qquad q^\al\;=\;-\,2\,\bn_\al\cdot\py\bH\:.
\ees
However, because $\bH=H^\al\bn_\al$\,, and the second fundamental form is symmetric, there holds
\begin{eqnarray*}
c-b&=&2\,\big(\bey\cdot\px-\bex\cdot\py\big)\,\bH\:\:=\:\:2\,H^\al\,\big(\bex\cdot\py-\bey\cdot\px\big)\,\bn_\al\\[1ex]
&=&2\,H^\al\,\text{e}^{\la}\,\big(h^\al_{21}-h^\al_{12}\big)\:\:=\:\:0\:.
\end{eqnarray*}
Whence, (\ref{a12}) may be recast in the form
\bes
\nabla^\perp\bL\;=\;\left(\begin{array}{cc}a&b\\[1ex]b&-\,a\end{array}\right)\left(\begin{array}{c}\bex\\[1ex]\bey\end{array}\right)\,-\,2\,\pro\big(\nabla\bH\big)\:.
\ees
Equivalently, in the $\,\{\bez,\bezz\}\,$ frame, this expression reads
\be\label{a0b}
\pz\bL\;=\;A\,\bezz\,-\,2\,i\,\pro\big(\pz\bH\big)\:,
\ee
where $\,A:=b+ia$\,.\\[1ex]
Using (\ref{a4}) and (\ref{a5}), the latter gives
\begin{eqnarray}\label{a1b}
\partial_{z^*z}\,\bL&=&\pzz\big(A\,\bezz\big)\,-\,2\,i\,\big(\bn_\al\cdot\pz\bH\big)\,\pzz\bn_\al\,-\,2\,i\,\pro\Big(\pzz\pro\big(\pz\bH\big)\Big)\nonumber\\[1ex]
&=&\Big[\text{e}^{-\la}\,\pzz\big(\text{e}^{\la}A\big)\,+\,2\,i\,\text{e}^{\la}\,\bH\cdot\pz\bH\Big]\,\bezz \,+\,2\,i\,\Big[\text{e}^{\la}\,\bH_0\cdot\pz\bH\Big]\,\bez\nonumber\\
&&\hspace{6cm}+\:\pro\,\partial_{z^*z}\bL\:.
\end{eqnarray}

\noindent
Observe that (\ref{a5}) and (\ref{a0b}) reveal that
\begin{eqnarray*}
&&\pro\,\partial_{z^*z}\,\bL\\[1ex]
&=&\Big[\pzz\big(\bn_\al\cdot\pz\bL\big)\,-\,\pzz\bn_\al\cdot\pz\bL\Big]\,\bn_\al\nonumber\\[1ex]
&=&\Big[-\,2\,i\,\pzz\Big(\bn_\al\cdot\pro\big(\pz\bH\big)\Big)\,+\,2\,i\,\pro\big(\pzz\bn_\al\big)\cdot\pro\big(\pz\bH\big)\,-\,\dfrac{1}{2}\,e^{\la}H_0^\al A\Big]\,\bn_\al\nonumber\\[1ex]
&=&2\,i\,\Big(\pro\big(\pzz\bn_\al\big)\,-\,\pzz\bn_\al\Big)\!\cdot\pro\big(\pz\bH\big)\,\bn_\al\,-\,2\,i\,\pro\,\pzz\pro\big(\pz\bH\big)\,-\,\dfrac{1}{2}\,e^{\la}\bH_0 A\nonumber\\[1ex]
&=&-\,2\,i\,\pro\,\pzz\pro\big(\pz\bH\big)\,-\,\dfrac{1}{2}\,e^{\la}\bH_0 A\:.
\end{eqnarray*}

\medskip
\noindent
Putting this into (\ref{a1b}) yields now
\begin{eqnarray}\label{a2b}
\partial_{z^*z}\,\bL&=&\Big[\text{e}^{-\la}\,\pzz\big(\text{e}^{\la}A\big)\,+\,2\,i\,\text{e}^{\la}\,\bH\cdot\pz\bH\Big]\,\bezz \,+\,2\,i\,\Big[\text{e}^{\la}\,\bH_0\cdot\pz\bH\Big]\,\bez\nonumber\\
&&\hspace{2cm}-\:2\,i\,\pro\,\pzz\pro\big(\pz\bH\big)\,-\,\dfrac{1}{2}\,e^{\la}\bH_0\, A\:.
\end{eqnarray}

\noindent
Because $\,\Im(\partial_{z^*z}\bL)=\vec{0}\,$, it suffices to ``project" the identity (\ref{a2b}) to discover 
\be\label{a3b}
\pzz\big(\text{e}^{\la}A\big)\:=\:-\,2\,i\,\text{e}^{2\la}\,\big(\bH\cdot\pz\bH\,+\,\bH^{*}_0\cdot\pzz\bH\big)
\ee
and
\be\label{a4b}
-\,4\,\Im\,\Big[i\,\pro\Big(\pzz\pro\big(\pz\bH\big)\Big)\Big]\;=\;e^{\la}\,\Im\big(A\,\vec{H}_0\big)\:.
\ee

\medskip
\noindent
Comparing (\ref{a3b}) to the Codazzi-Mainardi equation (\ref{a10}) shows that
\be\label{a5b}
A\;=\;-\,2\,i\,\text{e}^{\la}\,\bH_0^*\cdot\bH\,+\,i\,\text{e}^{-\la}f\:,
\ee
for some holomorphic function $f(z)$.\\[1.5ex]
\noindent
In conformal coordinates, there holds
\bes
\text{e}^{2\la}\,\Delta_\perp\bH\:=\:\pro\Big({div}\,\pro\big(\nabla\bH\big)\Big)\:\equiv\:-\,4\,\Im\Big[i\,\pro\Big(\pzz\pro\big(\pz\bH\big)\Big)\Big]\:,
\ees
so that (\ref{a4b}) becomes
\bes
\text{e}^{\la}\,\Delta_\perp\bH\;=\;\Im\big(A\,\vec{H}_0\big)\:.
\ees
Introducing (\ref{a5b}) in the latter gives
\be\label{a17}
\Delta_\perp\bH\,+\,2\,\Re\big(\bH\cdot\vec{H}^*_0\,\vec{H}_0\big)\;=\;\text{e}^{-2\la}\,\Im\big(i\,\vec{H}_0f\big)\:.
\ee
\medskip
To complete our derivation of the conformal Willmore equation, there remains to observe that\footnote{implicit summations over repeated indices are understood wherever appropriate.}
\begin{eqnarray*}
\Re\Big[H^\al_0\,\big(H^\bet_0\big)^{\!*}\Big]&=&\dfrac{1}{4}\,\Big[\big(h^\al_{11}-h^\al_{22}\big)\big(h^\bet_{11}-h^\bet_{22}\big)\,+\,4\,h^\al_{12}\,h^\bet_{12}\Big]\\[0ex]
&=&\dfrac{1}{2}\,\sum_{i,j=1}^{2}\,h^\al_{ij}\,h^\bet_{ij}\:-\:\dfrac{1}{4}\,\big(h^\al_{11}+h^\al_{22}\big)\big(h^\bet_{11}+h^\bet_{22}\big)\\[0ex]
&=&\dfrac{1}{2}\,\sum_{i,j=1}^{2}\,h^\al_{ij}\,h^\bet_{ij}\:-\:H^\al H^\bet\:.
\end{eqnarray*}
Accordingly, (\ref{a17}) becomes, after a few simple manipulations,
\bes
\Delta_\perp\bH\;+\,\sum_{i,j=1}^{2}\,h^\al_{ij}\,h^\bet_{ij}\,H^\bet\,\bn_\al\:-\:2\,\big|\vec{H}\big|^2\bH\:=\:\text{e}^{-2\la}\,\big\langle\vec{H}_0\,,f\big\rangle_{\mathbb{C}}\:.
\ees

\noindent
This identity is precisely the conformal Willmore equation.\\[-.5ex]

To finish the proof of Theorem~\ref{th-cw}, there remains to derive the identity (\ref{cwbis}). Let 
\be\label{trem0}
\nabla^\perp\bL_0\::=\:\nabla^\perp\bH\,-\,3\,\pi_{\bn}\big(\nabla\bH\big)\,+\,\star\,\big(\nabla^\perp\bn\wedge\bH\big)\:.
\ee
Since $\bH=H^\al\,\bn_\al$, we first note that
\begin{eqnarray}\label{trem1}
\nabla\bH\,-\,3\,\pi_{\bn}\big(\nabla\bH\big)&=&\big\langle\bex\,,\nabla\bH\big\rangle\,\bex\,-\,\big\langle\bey\,,\nabla\bH\big\rangle\,\bey\,-\,2\,\pi_{\bn}\big(\nabla\bH\big)\nonumber\\[1ex]
&=&H^\al\,\Big[\big\langle\bex\,,\nabla\bn_\al\rangle\,\bex\,+\,\big\langle\bey\,,\nabla\bn_\al\rangle\,\bey\Big]\,-\,2\,\pi_{\bn}\big(\nabla\bH\big)\nonumber\\[1ex]
&=&-\,H^\al\,\left(\begin{array}{cc}h^\al_{11}&h^\al_{21}\\[.5ex]h^\al_{12}&h^\al_{22}\end{array}\right)\nabla\bP\,-\,2\,\pi_{\bn}\big(\nabla\bH\big)\:.
\end{eqnarray}
It is shown in equation (II.33) from \cite{Riv2} that
\bes
\star\,\big(\nabla^\perp\bn\wedge\bH\big)\:=\:-\,H^\al\,\left(\begin{array}{cc}h^\al_{22}&-\,h^\al_{21}\\[.5ex]-\,h^\al_{12}&h^\al_{11}\end{array}\right)\nabla\bP\:.
\ees
Substituting the latter and (\ref{trem1}) into (\ref{trem0}) yields
\bes
\nabla^\perp\bL_0\::=\:-\,H^\al\,\left(\begin{array}{cc}h^\al_{11}-h^\al_{22}&2\,h^\al_{21}\\[.5ex]2\,h^\al_{12}&h^\al_{22}-h^\al_{11}\end{array}\right)\nabla\bP\,-\,2\,\pi_{\bn}\big(\nabla\bH\big)\:.
\ees
Equivalently, 
\bes
\pz\bL_0\:=\:-\,2\,i\,\text{e}^{\la}\,\bH^*_0\cdot\bH\;\bezz\,-\,2\,i\,\pi_{\bn}\big(\pz\bH\big)\:.
\ees
Introducing (\ref{a5b}), the latter may be recast as
\bes
\pz\bL_0\:=\:A\,\bezz\,-\,i\,\text{e}^{-\la}\,f(z)\,\bezz\,-\,2\,i\,\pi_{\bn}\big(\pz\bH\big)\:.
\ees
Accordingly, (\ref{a0b}) shows that
\bes
\pz\big(\bL_0-\bL\big)\:=\:-\,i\,\text{e}^{-\la}\,f(z)\,\bezz\:.
\ees
Finally, owing to (\ref{a4}) and to the holomorphy of $f$, we deduce the desired
\begin{eqnarray*}
\Delta\big(\bL_0-\bL\big)&\equiv&4\,\pzz\,\pz\big(\bL_0-\bL\big)\:\:=\:\:-\,i\,f(z)\,\pzz\big(e^{-\la}\,\bezz\big)\:\:=\:\:-\,2\,i\,\bH_0\,f\:.
\end{eqnarray*}

\reset
\appendix
\section{Appendix}

\subsection{Notational Conventions}

We append an arrow to all the elements belonging to $\R^m$. To simplify the notation, by $\bp\in X(\di)$ is meant $\bp\in X(\di,\R^m)$ whenever $X$ is a functional space. Similarly, we write $\nabla\bP\in \mathbb{R}^2\otimes X(\di)$ for $\nabla\bP\in X(\di,\mathbb{R}^{2m})\:.$\\[1.5ex]
Although this custom may seem odd at first, we allow the differential operators classically acting on scalars to act on elements of $\R^m$. Thus, for example, $\nabla\bp$ is the element of $\R^2\otimes\R^m$ that can be written $(\px\bp,\py\bp)$. If $S$ is a scalar and $\bR$ an element of $\R^m$, then we let
\begin{eqnarray*}
\bR\cdot\nabla\bP&:=&\big(\bR\cdot\px\bP\,,\,\bR\cdot\py\bP\big)\:,\\[1ex]
\nabla^\perp S\,\nabla\bP&:=&\px S\,\py\bp\,-\,\py S\,\px\bp\:,\\[1ex]
\nabla^\perp\bR\cdot\nabla\bP&:=&\px\bR\cdot\py\bp\,-\,\py\bR\cdot\px\bp\:,\\[1ex]
\nabla^\perp\bR\wedge\nabla\bP&:=&\px\bR\wedge\py\bp\,-\,\py\bR\wedge\px\bp\:.
\end{eqnarray*}
Similar quantities (with $\nabla$ in place of $\nabla^\perp$) are defined analogously, following the same logic. \\[1.5ex]
Two operations between multivectors are also useful. The interior multiplication $\res$ maps a pair comprising a $q$-vector $\gamma$ and a $p$-vector $\beta$ to a $(q-p)$-vector. It is defined via
\bes
\langle \gamma\res\beta\,,\alpha\rangle\;=\;\langle \gamma\,,\beta\wedge\alpha\rangle\:\qquad\text{for all $(q-p)$-vector $\al$.}
\ees
Let $\al$ be a $k$-vector. The first-order contraction operation $\bul$ is defined inductively through 
\bes
\al\bul\beta\;=\;\al\res\beta\:\:\qquad\text{when $\beta$ is a 1-vector}\:,
\ees
and
\bes
\al\bul(\beta\wedge\gamma)\;=\;(\al\bul\beta)\wedge\gamma\,+\,(-1)^{pq}\,(\al\bul\gamma)\wedge\beta\:,
\ees
when $\beta$ and $\gamma$ are respectively a $p$-vector and a $q$-vector. 

\subsection{Lorentz and Sobolev-Lorentz Spaces}

For the reader's convenience, we recall in this section the fundamentals of Lorentz spaces. More detailed accounts may be found in \cite{BL}, \cite {Hun}, and \cite{Ta1}.\\[-1ex]

For a real-valued measurable function $f$ on an open subset of $U\subset \mathbb{R}^n$, its belonging to a Lorentz space is determined by a condition involving the non-decreasing rearrangement of $|f|$ on the interval $(0,|U|)$, where $|U|$ denotes the Lebesgue measure of $U$. The non-increasing rearrangement $f^*$ of $|f|$ is the unique function from $(0,|U|)$ into $\mathbb{R}$ which is non-increasing and which satisfies
\bes
\Big|\big\{x\in U\;\big|\;|f(x)|\ge s\big\}\Big|\:=\:\Big|\big\{t\in(0,|U|)\;\big|\;f^*(t)\ge s\big\}\Big|\:.
\ees
\noindent
If $p\in(1,\infty)$ and $q\in[1,\infty]$, the Lorentz space $L^{p,q}(U)$ is the set of measurable functions $f:U\rightarrow\mathbb{R}$ for which
\bes
\int_{0}^{\infty}\,\big(t^{1/p}f^*(t)\big)^q\,\dfrac{dt}{t}\:<\:\infty\qquad\text{if}\:\:\:q<\infty\:,
\ees 
or
\bes
\sup_{t>0}\,t^{1/p}f^*(t)\:<\:\infty\qquad\text{if}\:\:\:q=\infty\:.
\ees \\
A norm on $L^{p,q}(U)$ is given by
\bes
\Vert f\Vert_{L^{p,q}(U)}\:=\:\big\Vert \,t^{1/p}f^{**}\big\Vert_{L^q([0,\infty),\,dt/t)}\qquad\text{where}\qquad f^{**}(t)\,:=\,\dfrac{1}{t}\,\int_{0}^{t}\,f^*(\tau)\,d\tau\:.
\ees
One verifies that
\bes
L^p\;=\;L^{p,p}\:,
\ees
and that $L^{p,\infty}$ is the weak-$L^p$ Marcinkiewicz space. \\[1ex]
Moreover, we have the inclusions
\bes
L^{p,1}\;\subset\;L^{p,q'}\;\subset\;L^{p,q''}\;\subset\;L^{p,\infty}\qquad\text{for}\qquad 1<q'<q''<\infty\:;
\ees
and if $U$ has finite measure, there holds for all $q$ and $q'$
\bes
L^{p',q'}(U)\;\subset\;L^{p,q}(U)\qquad\text{whenever}\qquad p<p'\:.
\ees
Finally, if $q<\infty$, the space $\,L^{\frac{p}{p-1},\frac{q}{q-1}}$\, is the dual of $L^{p,q}$.\\[1.5ex]
Similarly to Lebesgue spaces, Lorentz spaces obey a pointwise multiplication rule and a convolution product rule. More precisely, for $\,1<p_1,p_2<\infty\,$ and $\,1\le q_1,q_2\le\infty$, there holds
\bes
L^{p_1,q_1}\,\times\,L^{p_2,q_2}\:=\:L^{p,q}\qquad\text{with}\qquad\left\{\begin{array}{lcl}
p^{-1}&=&p_1^{-1}+\,p_2^{-1}\\[1ex]
q^{-1}&=&q_1^{-1}+\,q_2^{-1}\:,
\end{array}\right.
\ees
and
\bes
L^{p_1,q_1}\,*\,L^{p_2,q_2}\:=\:L^{p,q}\qquad\text{with}\qquad\left\{\begin{array}{lcl}
p^{-1}&=&p_1^{-1}+\,p_2^{-1}-\,1\\[1ex]
q^{-1}&=&q_1^{-1}+\,q_2^{-1}\:.
\end{array}\right.
\ees\\
\noindent
An interesting feature of Lorentz spaces is the possibility to generate them via interpolation of Lebesgue spaces. In particular, if $U\subset\mathbb{R}^n$ and $V\subset\mathbb{R}^m$ are open subsets, and if $r_0$, $r_1$, $p_0$, $p_1$ are real numbers satisfying
\bes
1\,\le\,r_0\,<\,r_1\,\le\,\infty\qquad\text{and}\qquad 1\,\le\,p_0\,\ne\,p_1\,\le\,\infty\:,
\ees
the following interpolation result holds. Let $T$ be a linear operator which\footnote{for $j\in\{0,1\}\:.$} maps continuously $L^{r_j}(U)$ into $L^{p_j}(V)$.
Then $T$ maps continuously $L^{r,q}(U)$ into $L^{p,q}(V)$ for each $\,q\in[1,\infty]$\, and every pair $(p,r)$ such that
\bes
\dfrac{1}{p}\;=\;\dfrac{1-\theta}{p_0}\,+\,\dfrac{\theta}{p_1}\qquad\text{and}\qquad\dfrac{1}{r}\;=\;\dfrac{1-\theta}{r_0}\,+\,\dfrac{\theta}{r_1}\:,
\ees
for some $\,\theta\in(0,1)$.\\[2ex]
\noindent
Lorentz spaces offer the possibility to sharpen the classical Sobolev embedding theorem. More precisely, it can be shown (see \cite{BL} and \cite{Ta3}) that
\bes
W^{k,q}(\mathbb{R}^m)\:\subset\:L^{p,q}(\mathbb{R}^m)
\ees
is a continuous embedding as long as
\bes
1\;\le\;q\;\le\;p\;<\;\infty\qquad\text{and}\qquad \dfrac{k}{m}\;=\;\dfrac{1}{q}\,-\,\dfrac{1}{p}\:.
\ees

\medskip
Our study requires that we introduce Sobolev-Lorentz spaces. These spaces are defined analogously to the ``standard" Sobolev spaces, but with the Lebesgue norms replaced by Lorentz norms. In an effort to simplify the presentation, we shall focus only on the two-dimensional unit disk $\di$. Let $m\in\mathbb{N}$, $p\in(1,\infty)$, and $q\in[1,\infty]$. The (homogeneous) Sobolev-Lorentz space $\,W^{m,(p,q)}(\di)\,$ consists of all locally summable functions $u$ on $\di$ such that $\,D^\al u\,$ exists in the weak sense and belongs to the Lorentz space $\,L^{p,q}(\di)\,$ for all multiindex $\,\al\,$ with $\,|\al|=m$. The norm
\bes
\Vert u\Vert_{W^{m,(p,q)}}\;:=\;\sum_{|\al|=m}\Vert D^\al u\Vert_{L^{p,q}}
\ees
clearly makes $\,W^{m,(p,q)}(\di)\,$ into a Banach space. The space $W^{0,(p,q)}$ is understood to be $L^{p,q}$. For notational convenience, we shall from now on omit to precise that we work on the unit disk $\di$. This shall arise no confusion. We also focus on a case of particular interest to us, namely $(p,q)=(2,1)$. 
Because $\,L^{2,1}$\, is a subspace of $\,L^2$\,, it follows immediately that for all $m\in\mathbb{N}$ there holds
\be\label{c0}
W^{m,(2,1)}\,\subset\,W^{m,2}\:,
\ee
where $W^{m,2}$ is the usual (homogeneous) Sobolev space.
In \cite{BW}, the authors prove\footnote{see also Theorem 3.3.4 in \cite{Hel}.} that $\,W^{1,(2,1)}$\, is a subspace of $\,L^\infty\cap\,C^0$, so that $\,W^{m,(2,1)}$\, is a subspace of $W^{m-1,\infty}\cap C^{m-1}$. Altogether, there holds
\be\label{c1}
W^{m,(2,1)}\,\subset\,W^{m,2}\cap W^{m-1,\infty}\cap C^{m-1}\:.
\ee
The Rellich-Kondrachov theorem states that $\,W^{1,2}\,$ is compact in $\,L^r\,$ for all finite $\,r\ge1\,.$ By standard interpolation techniques, it ensues that
\be\label{c11}
W^{m,2}\,\subset\subset\,W^{m-1,(p,q)}\qquad\text{for}\qquad\left\{\begin{array}{l}p=1=q\\[1ex]
1<p<\infty\:\:,\:\:1\le q\le\infty\:.\end{array}\right.
\ee
are compact inclusions for every $m\in\mathbb{N}^*$. \\[2ex]
The embedding (\ref{c0}) shows that each element of $\,W^{m,(2,1)}\,$ has a well-defined trace on the boundary of the unit disk, for $\,m\ge1$. If that trace is null, the said element belongs to the space $\,W_0^{m,(2,1)}$. In addition, it is easily seen that $W_0^{m,(2,1)}$ is the closure of $C_c^{\infty}$ in $W^{m,(2,1)}$.\\[1ex]
It is instructive to characterize the dual of the space $W_0^{m,(p,q)}$, which we shall denote $\,W^{-m,(p',q')}$, where $(p',q')$ is the conjugate pair of $(p,q)$\,:
\be\label{y4}
p'\;=\;\big(1\,-\,p^{-1}\big)^{-1}\qquad\text{and}\qquad q'\;=\;\big(1\,-\,q^{-1}\big)^{-1}\:.
\ee

\begin{Prop}
\label{charsob}
Suppose that $T$ is an element of $\,W^{-m,(p',q')}(\di,\mathbb{R})\,,$ so that $T$ is a bounded linear functional on $\,W_0^{m,(p,q)}$. For every $\delta>0$, there exists an element $P$ of $\,W^{1-m,(p',q')}(\di,\mathbb{R}^2)\,$ such that
\be\label{c2}
T\;=\;div\;P\:,
\ee
and
\be\label{c3}
\big\Vert P\big\Vert_{W^{1-m,(p',q')}}\;\le\;\Vert T\Vert_{W^{-m,(p',q')}}\,+\,\delta\:.
\ee
\end{Prop}
$\textbf{Proof.}$ This is a mere adaption of the analogous statement for the Lebesgue-Sobolev space $W^{m,p}$, ultimately following from the Hahn-Banach theorem. The reader is referred to Theorem 3.10 in \cite{Ad1} for details.\\[-3ex]

$\hfill\blacksquare$ \\

We next bring into light an interesting Hodge decomposition result which follows from standard elliptic theory and the interpolation nature of Lorentz spaces (cf. Proposition 3.3.9 in \cite{Hel}). 

\begin{Prop}
\label{prop-hel}
To every vector field $\,\vec{g}=(g_1,g_2)\in L^1(\di,\mathbb{R}^2)\,$ we associate the functions $\al$, $\beta$, and $h$ on $\di$ which are solutions of
\bes
\left\{\begin{array}{rcl}
\Delta\al&=&div\,\vec{g}\\[1ex]
\al&=&0\end{array}\right.\quad,\quad
\left.\begin{array}{rclcl}
\Delta\beta&=&curl\,\vec{g}&\quad&\text{in\:\:\:}\di\\[1ex]
\beta&=&0&\quad&\text{on\:\:\:}\partial\di\:,
\end{array}\right.
\ees
and
\bes
\left\{\begin{array}{rcl}
g_1&=&\dfrac{\partial\al}{\partial x}\,+\,\dfrac{\partial\beta}{\partial y}\,+\,\dfrac{\partial h}{\partial x}\\[3ex]
g_2&=&\dfrac{\partial\al}{\partial y}\,-\,\dfrac{\partial\beta}{\partial x}\,+\,\dfrac{\partial h}{\partial y}\:,
\end{array}\right.
\ees
so that $h$ is harmonic on $\di$. \\
Then the operators
\bes
\vec{g}\;\longmapsto\;\nabla\al\:,\qquad \vec{g}\;\longmapsto\;\nabla\beta\:,\qquad\text{and}\qquad \vec{g}\;\longmapsto\;\nabla h
\ees
map continuously $L^{p,q}(\di,\mathbb{R}^2)$ into itself for every choice of $p\in(1,\infty)$ and $q\in[1,\infty]$. 
\end{Prop}

\noindent
This lemma is the central ingredient to obtain the following result. 

\begin{Lma}
\label{A-1}
Let $p\in(1,\infty)$ and $q\in(1,\infty]$. Suppose that $G$ is an element of $W^{-1,(p,q)}(\di,\mathbb{R}^2)$ which satisfies (in the functional sense)
\bes
div\,G\;=\;0\qquad\text{in}\quad\di\:.
\ees
Then there exists an element $L$ in the space $L^{p,q}(\di,\mathbb{R})$ such that
\be\label{y1}
G\;=\;\nabla^{\perp}L\:.
\ee
\end{Lma}
$\textbf{Proof.}$ Let $p'$ and $q'$ be as in (\ref{y4}). We remind the reader that $W^{-1,(p,q)}$ is the dual of $W_0^{1,(p',q')}$.\\
A classical result of Laurent Schwartz guarantees the existence of $L$ in $\mathcal{D}'$ such that (\ref{y1}) holds in the sense of distributions.  
Owing to the density of $\,\mathcal{D}=C^\infty_c$ in $W_0^{1,(p',q')}$, it thus follows that
\be\label{y2}
\int_{\di} G\cdot F\:=\:-\,\int_{\di}L\;div\,F\:,\qquad\forall\:\: F\,\in\,W_0^{1,(p',q')}(\di,\mathbb{R}^2)\:.
\ee
Next, let $f$ be an arbitrary element of $L^{p',q'}(\mathbb{R})$ with null average over the unit-disk. 
Choosing $\,\vec{g}=(f,0)\,$ in Proposition~\ref{prop-hel}, we infer the existence of an element $F:=(\al,\beta)\,$ in $W_0^{1,(p',q')}(\di,\mathbb{R}^2)$ satisfying
\be\label{y3}
div\,F\;=\;f\qquad\text{and}\qquad \Vert F\Vert_{W^{1,(p',q')}}\:\le\:\Vert f \Vert_{L^{p',q'}}\:.
\ee
Altogether, (\ref{y2}) and (\ref{y3}) show that $L$ acts linearly on $L^{p',q'}$ with the estimate
\bes
\bigg|\int_{\di}L\,f\,\bigg|\;\;\le\;\;\Vert G\Vert_{W^{-1,(p,q)}}\,\Vert f\Vert_{L^{p',q'}}\:.
\ees
Since $L^{p,q}$ is the dual space of $L^{p',q'}$, the desired statement ensues. \\[-3ex]

$\hfill\blacksquare$

\subsection{Miscellaneous Identities}

\begin{Lma}
\label{A-2}
Let $\bp$, $\bn$, and $\bH$ be as in Theorem~\ref{th-ps}. We set
\bes
Q\;:=\;\nabla\bH\,-\,3\,\pro\big(\nabla\bH\big)\,+\,\star\big(\nabla^\perp\bn\wedge\bH\big)\:.
\ees
Then the following identities hold
\be\label{b1}
\nabla\bp\cdot Q\;=\;0\hspace{1cm}\text{and}\hspace{1cm}\nabla\bp\wedge Q\;=\;-\,2\,\nabla\bp\wedge\nabla\bH\:.
\ee
\end{Lma}
\smallskip
$\textbf{Proof.}$
Firstly, we note that
\bes
\nabla^\perp\bna\;=\;\big\langle\bex\,,\nabla^\perp\bna\big\rangle\,\bex\,+\,\big\langle\bey\,,\nabla^\perp\bna\big\rangle\,\bey\,+\,\big\langle\bn_\beta\,,\nabla^\perp\bna\big\rangle\,\bn_\beta\:,
\ees
so that
\bes
\star\big(\bn\wedge\nabla^\perp\bna\big)\;=\;\big\langle\bex\,,\nabla^\perp\bna\big\rangle\,\bey\,-\,\big\langle\bey\,,\nabla^\perp\bna\big\rangle\,\bex\:.
\ees
Whence,
\bes
\star\big(\nabla^\perp\bn\wedge\bna\big)\;=\;\left(\begin{array}{c}-\,h^\al_{22}\\[1ex]h^\al_{12}\end{array}\right)\px\bp\;+\,\left(\begin{array}{c}h^\al_{12}\\[1ex]-\,h^\al_{11}\end{array}\right)\py\bp\:.
\ees
Accordingly, we find
\bes
\nabla\bp\wedge\,\star\big(\nabla^\perp\bn\wedge\bna\big)\;=\;\big(h^\al_{12}\,-\,h^\al_{12}\big)\,\px\bp\wedge\py\bp\;=\;0\:,
\ees
and
\bes
\nabla\bp\,\cdot\,\star\big(\nabla^\perp\bn\wedge\bna\big)\;=\;-\,\text{e}^{2\la}\,\big(h^\al_{11}\,+\,h^\al_{22}\big)\;=\;-\,2\,\text{e}^{2\la}\,H^\al\:.
\ees\\[-1ex]
The last two identities yield
\be\label{b2}
\nabla\bp\wedge\,\star\big(\nabla^\perp\bn\wedge\bH\big)\;=\;0\:,
\ee
and
\be\label{b3}
\nabla\bp\,\cdot\,\star\big(\nabla^\perp\bn\wedge\bH\big)\;=\;-\,2\,\text{e}^{2\la}\,\big|\bH\big|^2\:.
\ee\\
\noindent
For three indices $\,(a,b,c)\in\{x^1,x^2\}\,$, let
\bes
F(a,b,c)\;:=\;\partial_a\bp\,\wedge\,\big\langle\beb\,,\partial_c\bH\big\rangle\,\beb\:.
\ees
We have
\begin{eqnarray*}
F(a,b,c)&=&e^\la\,\big\langle\beb\,,\partial_c\bH\big\rangle\,\bea\wedge\beb\nonumber\\[1ex]
&=&-\,\text{e}^\la\,H^\al\,\big\langle\bna\,,\partial_c\beb\big\rangle\,\bea\wedge\beb\hspace{1cm}\text{since \,$\bH=H^\al\bna\:$}\nonumber\\[1ex]
&=&-\,\text{e}^\la\,H^\al\,h^\al_{cb}\,\bea\wedge\beb\:.
\end{eqnarray*}
Using this, we obtain
\begin{eqnarray*}
\nabla\bp\,\wedge\,\Big(\nabla\bH\,-\,\pro\big(\nabla\bH\big)\Big)&\equiv&\nabla\bp\,\wedge\,\Big(\big\langle\bex\,,\nabla\bH\big\rangle\,\bex\,+\,\big\langle\bey\,,\nabla\bH\big\rangle\,\bey\Big)\nonumber\\[1ex]
&=&F(x^1,x^1,x^1)\,+\,F(x^1,x^2,x^1)\nonumber\\[.5ex]
&&\hspace{1cm}+\:F(x^2,x^1,x^2)\,+\,F(x^2,x^2,x^2)\nonumber\\[1ex]
&=&-\,\text{e}^\la\,H^\al\,h^\al_{21}\,\bex\wedge\bey\,-\,\text{e}^\la\,H^\al\,h^\al_{12}\,\bey\wedge\bex\nonumber\\[1ex]
&=&0\:.
\end{eqnarray*}
Hence
\be\label{b4}
\nabla\bp\,\wedge\,\Big(\nabla\bH\,-\,3\,\pro\big(\nabla\bH\big)\Big)\;=\;-\,2\,\nabla\bp\,\wedge\,\nabla\bH\:.
\ee\\[-1ex]
An evident yet nonetheless useful identity is
\bes
\nabla\bp\,\cdot\pro\big(\nabla\bH\big)\;=\;0\:,
\ees
thereby giving
\be\label{b5}
\nabla\bp\,\cdot\,\Big(\nabla\bH\,-\,3\,\pro\big(\nabla\bH\big)\Big)\;=\;\nabla\bp\,\cdot\,\nabla\bH\:.
\ee\\
Bringing altogether (\ref{b2}) and (\ref{b4}) yields the second part of (\ref{b1}). Deriving the first part of (\ref{b1}) requires a bit more work. Combining (\ref{b3}) and (\ref{b5}) shows that
\be\label{b6}
\nabla\bp\,\cdot\,Q\;=\;\nabla\bp\,\cdot\,\nabla\bH\,-\,2\,\text{e}^{2\la}\,\big|\bH\big|^2\:.
\ee
For three indices $\,(a,b,c)\in\{x^1,x^2\}\,$, let
\bes
G(a,b,c)\;:=\;\partial_a\bp\,\cdot\,\big\langle\beb\,,\partial_c\bH\big\rangle\,\beb\:.
\ees
Just as above, we verify easily that
\bes
G(a,b,c)\;=\;\text{e}^{2\la}\,H^\al\,h^\al_{cb}\,\delta_{ab}\:.
\ees
Consequently, there holds
\begin{eqnarray*}
\nabla\bp\,\cdot\,\nabla\bH&=&G(x^1,x^1,x^1)\,+\,G(x^1,x^2,x^1)\,+\,G(x^2,x^1,x^2)\,+\,G(x^2,x^2,x^2)\\[1ex]
&=&\text{e}^{2\la}\,H^\al\,h^\al_{11}\,+\,\text{e}^{2\la}\,H^\al\,h^\al_{22}\\[1ex]
&=&2\,\text{e}^{2\la}\,\big|\bH\big|^2\:.
\end{eqnarray*}
Putting this into (\ref{b6}) finally gives the sought after first identity in (\ref{b1}).\\[-1.5ex]

$\hfill\blacksquare$

\begin{Lma}
\label{A-3} 
Using the notation introduced in Section III.2.2, the Codazzi-Mainardi identity may be recast in the form
\be\label{a10}
\text{e}^{-2\la}\,\pzz\big(\text{e}^{2\la}\,\bH^*_0\cdot\bH\big)\;=\;\bH\cdot\pz\bH\,+\,\bH^*_0\cdot\pzz\bH\:.
\ee
\end{Lma}
$\textbf{Proof.}$ Consider the 1-form
\bes
\eta_\al\;:=\;\big\langle\bezz, d\bn_\al\big\rangle\:.
\ees
On one hand, there holds
\begin{eqnarray*}
\eta_\al&=&\big(\bezz\cdot\pz\bn_\al\big)\,dz\,+\,\big(\bezz\cdot\pzz\bn_\al\big)\,dz^*\:\:\stackrel{(\ref{a5})}{=}\:\:-\,\dfrac{\text{e}^\la}{2}\,\big(H^\al dz+H^\al_0dz^*\big)\:,
\end{eqnarray*}
so that
\be\label{a6}
d\eta_\al\;=\;-\,\dfrac{\text{e}^\la}{2}\,\big(H^\al_0\pz\la\,-\,H^\al\pzz\la\,+\,\pz H^\al_0\,-\,\pzz H^\al\big)\,dz\wedge dz^*\:.
\ee
On the other hand, with (\ref{a4}) and (\ref{a5}), we find
\begin{eqnarray*}
d\eta_\al&=&d\,\big\langle\bezz ,d\bn_\al\big\rangle\:\:=\:\:\big(\pz\bezz\cdot\pzz\bn_\al\,-\,\pzz\bezz\cdot\pz\bn_\al\big)\,dz\wedge dz^*\nonumber\\[1ex]
&=&\dfrac{\text{e}^\la}{2}\,\big(H^\al_0\pz\la\,+\,H^\al\pzz\la\,+\,\vec{H}\cdot\pzz\bn_\al\,-\,\vec{H}_0\cdot\pz\bn_\al\big)\,dz\wedge dz^*\:.
\end{eqnarray*}
Comparing the latter to (\ref{a6}) yields
\be\label{a7}
2\,H_0^\al\,\partial_z\la\;=\;\pzz H^\al\,-\,\pz H^\al_0\,-\,\vec{H}\cdot\pzz\bn_\al\,+\,\vec{H}_0\cdot\pz\bn_\al\:.
\ee
By antisymmetry, there holds\footnote{we ``contextually" sum over repeated indices whenever appropriate.}
\bes
H^\al\vec{H}\cdot\pzz\bn_\al\;=\;H^\al H^\beta\,\bn_\beta\cdot\pzz\bn_\al\;=\;0\:,
\ees
and whence
\be\label{a8}
H^\al\pzz H^\al\;=\;H^\al\pzz H^\al\,+\,H^\al\vec{H}\cdot\pzz\bn_\al\;=\;\vec{H}\cdot\pzz\vec{H}\:.
\ee
Similarly, we find
\begin{eqnarray}\label{a9}
H^\al\big(\pz H^\al_0\,-\,\vec{H}_0\cdot\pz\bn_\al\big)&=&\pz\big(\vec{H}_0\cdot\vec{H}\big)\,-\,H_0^\al\pz H^\al\,-H^\al\vec{H}_0\cdot\pz\bn_\al\nonumber\\[1ex]
&=&\pz\big(\vec{H}_0\cdot\vec{H}\big)\,-\,\vec{H}_0\cdot\pz\vec{H}\:.
\end{eqnarray}
Multiplying (\ref{a7}) throughout by $H^\al$, summing over $\al$, and using (\ref{a8}) and (\ref{a9}) gives after a few manipulations the (complex conjugate of the) desired
\bes
\text{e}^{-2\la}\,\pz\big(\text{e}^{2\la}\,\vec{H}_0\cdot\vec{H}\big)\;=\;\vec{H}\cdot\pzz\vec{H}\,+\,\vec{H}_0\cdot\pz\vec{H}\:.
\ees\\[-8ex]

$\hfill\blacksquare$

\begin{Lma}
\label{A-4} Using the notation of Section III.2, there holds
\be\label{r1b}
-\,\Delta\bp\;=\;\nabla\bR\,\bul\nabla^\perp\bp\;+\,\nabla S\;\nabla^\perp\bp\:.
\ee
\end{Lma}
$\textbf{Proof.}$ Note that for any 1-vector $\ba$, we have
\bes
(\ba\wedge\bej)\bul\bei\;=\;(\bei\res\ba)\wedge\bej\,+\,\ba\wedge(\bei\res\bej)\;=\;(\bei\cdot\ba)\,\bej\,+\,\delta_{ij}\,\ba\:.
\ees
From this, and $\,\bei:=\text{e}^{-\la}\partial_{x^i}\bp\,$, it follows easily that whenever
\bes
\bV\;:=\;V^{i}\,\bei\,+\,V^{\al}\,\bna\:,
\ees\\[-5ex]
then
\bes
\left\{\begin{array}{rcccl}
\big(\bV\wedge\nabla^\perp\bp\big)\bul\nabla^\perp\bp&=&\text{e}^{2\la}\,\big(3\,V^i\,\bei\,+\,2\,V^\al\,\bna\big)&&\\[1.5ex]
\big(\bV\wedge\nabla\bp\big)\bul\nabla^\perp\bp&=&\text{e}^{2\la}\,\big(V^{2}\,\bex\,-\,V^{1}\,\bey\big)&\equiv&\big(\bV\cdot\nabla\bp\big)\nabla^\perp\bp\:.
\end{array}\right.
\ees
In particular, since $\,\bH=H^\al\,\bna\,$, we find that
\begin{eqnarray*}
\nabla\bR\,\bul\nabla^\perp\bp&\equiv&-\,\big(\bL\wedge\nabla\bp\,+\,2\,\bH\wedge\nabla^\perp\bp\big)\,\bul\nabla^\perp\bp\\[1ex]
&=&-\,\big(\bL\cdot\nabla\bp\big)\nabla^\perp\bp\,-\,2\,\text{e}^{2\la}\,\bH\\[1ex]
&=&-\,\nabla S\;\nabla^\perp\bp\,-\,2\,\text{e}^{2\la}\,\bH\:.
\end{eqnarray*}
Hence,
\bes
\nabla\bR\,\bul\nabla^\perp\bp\;+\,\nabla S\;\nabla^\perp\bp\;=\;-\,2\,\text{e}^{2\la}\,\bH\:.
\ees

\medskip
\noindent
Finally, there remains to observe that
\bes
\Delta\bp\;=\;2\,\text{e}^{2\la}\,\bH
\ees
to deduce the desired identity (\ref{r1b})\:.\\[-4ex]

$\hfill\blacksquare$

\eject

\end{document}